\documentclass[10pt, dvipsnames]{amsart}
\usepackage{amsmath,amssymb,amsfonts, mathtools}
\usepackage[mathscr]{eucal}
\usepackage[margin=1in]{geometry}
\usepackage[all]{xy}
\usepackage{bbm}
\usepackage{xparse}

\usepackage{tikz, tikz-cd}%
\usetikzlibrary{matrix,arrows,decorations.pathmorphing,cd,patterns,calc,backgrounds,patterns}
\tikzset{dummy/.style= {circle,fill,draw,inner sep=0pt,minimum size=1.2mm}}
\tikzset{vertex/.style={fill, circle, minimum size=.1cm, inner sep=0pt}}

\usepackage[pagebackref, colorlinks, citecolor=DarkOrchid, linkcolor=black, urlcolor=NavyBlue]{hyperref}
\hypersetup{linktoc=all,} 

\hypersetup{
  colorlinks   = true, 
  urlcolor     = blue, 
  linkcolor    = PineGreen, 
colorlinks, citecolor    = DarkOrchid 
}

\usepackage{comment}

\usepackage[textwidth=25mm, textsize=tiny]{todonotes}
\usepackage[final]{microtype}

\usepackage{enumerate}

\usepackage{appendix}

\numberwithin{equation}{section} 
\numberwithin{figure}{section}
\usepackage[nameinlink,capitalise,noabbrev]{cleveref}

\newcommand{\newrefformat}[2]{}

\usepackage{stmaryrd}

\newcommand\restr[2]{{
  \left.\kern-\nulldelimiterspace 
  #1 
  \vphantom{\big|} 
  \right|_{#2} 
  }}

\crefname{lemma}{Lemma}{Lemmas}
\crefname{theorem}{Theorem}{Theorems}
\crefname{definition}{Definition}{Definitions}
\crefname{proposition}{Proposition}{Propositions}
\crefname{remark}{Remark}{Remarks}
\crefname{corollary}{Corollary}{Corollaries}
\crefname{question}{Question}{Questions}
\crefname{equation}{Equation}{Equations}
\crefname{construction}{Construction}{Constructions}
\crefname{ex}{Example}{Examples}
\crefname{appsec}{Appendix}{Appendices}
\crefname{subsection}{Subsection}{Subsections}
\crefname{section}{Section}{Sections}

\theoremstyle{plain}
\newtheorem{theorem}[equation]{Theorem}
\newtheorem{corollary}[equation]{Corollary}
\newtheorem{proposition}[equation]{Proposition}

\theoremstyle{definition}
\newtheorem{definition}[equation]{Definition}
\newtheorem{example}[equation]{Example}

\newtheorem{construction}[equation]{Construction}
\newtheorem{remark}[equation]{Remark}

\newtheorem{introtheorem}{Theorem}
 
\crefname{introtheorem}{Theorem}{Theorems}

\usepackage[textwidth=25mm, textsize=tiny]{todonotes}

\newcommand{\coHH}{\mathrm{coHH}}
\newcommand{\THH}{\mathrm{THH}}
\newcommand{\coTHH}{\mathrm{coTHH}}

\newcommand{\HH}{\mathrm{HH}}
\newcommand{\Hom}{\mathrm{Hom}}
\newcommand{\coMod}{\mathrm{coMod}}
\newcommand{\Mod}{\mathrm{Mod}}

\newcommand{\Vect}{\mathrm{Vect}}
\newcommand{\op}{\mathrm{op}}

\newcommand{\coAlg}{\mathrm{coAlg}}
\newcommand{\Alg}{\mathrm{Alg}}

\newcommand{\cC}{\mathscr{C}}
\newcommand{\cD}{\mathscr{D}}

\newcommand{\dkx}{k\langle X\rangle}
\newcommand{\dky}{k\langle \langle Y \rangle \rangle}

\newcommand{\FF}{k}

\newcommand{\Lin}{\mathrm{Lin}}

\newcommand{\Proj}{\mathrm{Proj}}
\newcommand{\Inj}{\mathrm{Inj}}

\newcommand{\tr}{\mathrm{tr}}

\newcommand{\cotr}{\mathrm{cotr}}

\newcommand{\fd}{\mathrm{fd}}
\newcommand{\fc}{\mathrm{fc}}
\newcommand{\fg}{\mathrm{fg}}
\newcommand{\lf}{\mathrm{lf}}

\newcommand{\SI}{\Sigma^\infty_+}
\newcommand{\id}{\mathrm{id}}

\newcommand{\Map}{\mathrm{Map}}

\newcommand{\GL}{\mathrm{GL}}

\newcommand{\cK}{K^c}
\newcommand{\bK}{G^c}

\NewDocumentCommand{\cotens}{e{_^}}{%
  \mathbin{\mathop{\square}\displaylimits
    \IfValueT{#1}{_{#1}}
    \IfValueT{#2}{^{#2}}
  }%
}

\NewDocumentCommand{\tens}{e{_^}}{%
  \mathbin{\mathop{\otimes}\displaylimits
    \IfValueT{#1}{_{#1}}
    \IfValueT{#2}{^{#2}}
  }%
}


\begin{document}

\author[T. Gerhardt]{Teena Gerhardt}
\address{Department of Mathematics, Michigan State University, 619 Red Cedar Rd, East Lansing, MI 48824, USA}
\email{teena@math.msu.edu}

\author[M.\ P\'eroux]{Maximilien P\'eroux}
\address{Department of Mathematics, Michigan State University, 619 Red Cedar Rd, East Lansing, MI 48824, USA}
\email{peroux@msu.edu}

\author[W.H.B. Sor\'e]{W. Hermann B. Sor\'e}
\address{D\'epartement de Math\'ematiques et Informatique, Universit\'e Nazi Boni de Bobo-Dioulasso, 01 BP 1091 Bobo-Dioulasso 01, Burkina Faso}
\email{hermann.sore@gmail.com}

\title[Coalgebraic $K$-theory]{Coalgebraic $K$-theory}

\begin{abstract}
    We establish comparison maps between the classical algebraic $K$-theory of algebras over a field and its analogue $K^c$, an algebraic $K$-theory for coalgebras over a field.
    The comparison maps are compatible with the Hattori--Stallings (co)traces.
    We identify conditions on the algebras or coalgebras under which the comparison maps are equivalences. 
    Notably, the algebraic $K$-theory of the power series ring is equivalent to the $K^c$-theory of the divided power coalgebra.
    We also establish comparison maps between the $G$-theory of finite dimensional representations of an algebra and its analogue $G^c$ for coalgebras. In particular, we show that the Swan theory of a group is equivalent to the $G^c$-theory of the representative functions coalgebra, reframing the classical character of a group as a trace in coHochschild homology.
\end{abstract}

\maketitle

\section{Introduction}
Algebraic $K$-theory is a foundational invariant for rings and its study sheds light on various areas of mathematics, including algebraic topology, algebraic geometry, and number theory.  One fruitful approach to the study of algebraic $K$-theory is trace methods, in which algebraic $K$-theory is approximated by invariants which are more computable. For a ring $R$, there is a map from its algebraic $K$-theory to its Hochschild homology, called the Dennis trace map $K_n(R)\to \HH_n(R)$. Goodwillie \cite{TraceGood} proved that the Dennis trace factors through  
negative cyclic homology, $\mathrm{HC}^{-}_n(R)$, 
which is often a good approximation of algebraic $K$-theory rationally.
To acquire integral information about algebraic $K$-theory, topological versions of Hochschild homology and cyclic homology were developed \cite{bok}, \cite{BHM}, as well as a trace map from algebraic $K$-theory to topological Hochschild homology (THH), which factors through topological cyclic homology (TC). Topological cyclic homology often offers a very good approximation to $K(R)$ (see, for example, \cite{DGM}). This trace method approach has been hugely successful in studying algebraic $K$-theory.

Coalgebras, the dual to rings and algebras, appear throughout algebraic topology. For instance, the Landweber exactness theorem implies that complex oriented homologies admit natural comodule structures over a coalgebra determined by the complex bordism spectrum  $MU$ \cite{Landweber}. Similarly, by the Milnor--Moore theorem, the rational homology of a based loop space and the universal enveloping algebra of its rational homotopy Lie algebra are isomorphic as Hopf algebras \cite{MilnorMoore}.
Further, many famous dualities in topology capture interplay between algebraic and coalgebraic structures.

In the setting of trace theories, such interesting dualities arise. For coalgebras over a field $k$ there is a theory of coHochschild homology (coHH), dual to Hochschild homology, due to Doi \cite{Doi}.
Hess and Shipley extended the definition of coHH into the topological setting, and defined topological coHochschild homology, $\coTHH(C)$, for any coalgebraic spectrum $C$ \cite{HScothh}, see also \cite{cothhstring, dualitySW}. 
Dualities between algebras and coalgebras lead to relationships between $\THH$ and $\coTHH$, providing new approaches to the computation of $\THH$.
For example, Hess--Shipley showed that given some conditions on a connected space $X$, there is an equivalence \[\THH(\SI \Omega X)\simeq \SI \mathscr{L}X\simeq \coTHH(\SI X),\] when viewing the suspension spectrum $\SI X$ as a coalgebra spectrum via the diagonal $X\to X\times X$ \cite{HScothh}. This provides a new model for the free loop space $\mathscr{L}X$ of $X$. 
Using a dual version of the B\"okstedt 
spectral sequence on $\coTHH$ \cite{toolscothh}, the  first author, together with Bohmann and Shipley provided new computations of the homologies of $\mathscr{L}X$ \cite{cothhstring} that were less accessible using classical $\THH$.
The second author with Bay{\i}nd{\i}r established in \cite{dualitySW} another duality comparison between the coTHH of a coalgebra spectrum $C$ and the THH of the Spanier--Whitehead dual of $C$.

As these dualities give new insights on $\THH$, this naturally leads to the question: is there an analogue of algebraic $K$-theory for coalgebras? Such a theory could lead to new insights on the usual algebraic $K$-theory of rings via duality.
To address this question, together with Klanderman, the second author introduced the \textit{coalgebraic} $K$-theory $K^c(C)$ of a coalgebra $C$ over a field $k$ in \cite{cotrace}. This $K$-theory is defined as the algebraic $K$-theory of the exact category of finitely cogenerated and injective comodules over $C$ (see also \cref{section: intro of K^c and G^c} below). 
The linear dual $C^*=\Hom_k(C, k)$ is naturally a $k$-algebra, and in this paper, we establish a comparison between the algebraic $K$-theory of the algebra $C^*$ and the coalgebraic $K$-theory $K^c(C)$.

\begin{introtheorem}[\cref{thm coalg KT to alg KT}, \cref{thm: comparison map of K^c is ring hm}, \cref{corollary: equivalence  K(C)=K(C^*)}]
Let $k$ be a field.
   Given a $k$-coalgebra $C$, there is a natural map of algebraic $K$-theory spectra:
\[
(-)^*\colon \cK(C)\longrightarrow K(C^*).
\]
 If $C$ is cocommutative, this map is a ring homomorphism of commutative ring spectra. Further, if $C$ is finite dimensional, or if $C^*$ is a commutative Noetherian $k$-algebra, this map is an equivalence.
\end{introtheorem}

There is an analogous result given a $k$-algebra $A$ rather than a coalgebra, however, $A^*$ might not be a $k$-coalgebra in general. Instead we consider $A^\circ\subseteq A^*$, the largest subvector space so that the dual of the multiplication becomes a comultiplication. This is often referred as the finite dual or Sweedler dual of $A$ \cite{Sweedler}. In this case, there is another comparison map $K(A)\to K^c(A^\circ)$ when $A$ is commutative and Noetherian, see \cref{theorem: finite dual on algebraic K theory of ring}, that becomes an equivalence for certain algebras, see \cref{theorem: equivalence K(A)=K(A^*)}.

Most notably, a particular example of interest is the $k$-algebra given by the power series $k[[y]]$. 
Its finite dual $k[[y]]^\circ$ is precisely the divided power coalgebra $\dkx$, and in characteristic zero $\dkx \cong k[x]$ (see \cref{example: k[x] is coreflexive}). 
Here $X$ denotes the linear transformation $k[[y]] \to k$ that projects to the coefficient of $y$.
Therefore, a consequence of the results above is the following.

\begin{introtheorem}[\cref{corollary: K-theory of power series}]
Let $k$ be a field.
    There is a natural equivalence of algebraic $K$-theory ring spectra:
    \[
    \cK(\dkx)\simeq K(k[[y]]).
    \]
\end{introtheorem}

It is our hope that this coalgebraic framework provides new insights into the algebraic $K$-theory of power series, where the uncountable dimensional algebra $k[[y]]$ is replaced with a countable dimensional coalgebra $\dkx$.
This  hope motivates the need for computationally effective tools to determine the coalgebraic $K$-theory of coalgebras. In particular, one would like to develop trace methods for coalgebraic $K$-theory. Klanderman and the second author established an analogue of the Dennis trace $K_n(R)\to \HH_n(R)$ in the coalgebraic setting for $n=0$ \cite{cotrace}. 
In particular, they introduced a coalgebraic analogue of the Hattori--Stallings trace $K^c_0(C)\to (\coHH_0(C))^*$, for any $k$-coalgebra $C$.
Our next result shows that the Hattori--Stallings traces in the algebraic and coalgebraic settings are compatible.

\begin{introtheorem}[\cref{theorem: cotrace compatible with K(C^*)}, \cref{theorem: finite dual on algebraic K theory of ring}]\label{theorem:comparison}
Let $k$ be a field.
    The comparison maps $K^c(C)\to K(C^*)$ and $K(A)\to K^c(A^\circ)$ established above for $k$-coalgebras $C$ and $k$-algebras $A$ are compatible with the Hattori--Stallings (co)traces. 
\end{introtheorem}

The result above also suggests a definition for trace maps $K_n^c(C)\to \coHH_n(C)^*$ for all $n\geq 1$. Forthcoming work of the second author with Agarwal and Mehrle aims to further study  these trace maps (see \cref{Rem: higher trace maps}).

The use of a coalgebraic framework in algebraic $K$-theory has appeared before in the literature, in a different context. 
Recall that the stable geometry of a space $X$ is captured by Waldhausen's $A$-theory, denoted $A(X)$, defined as the algebraic $K$-theory of the Waldhausen category formed by finitely dominated retractive spaces over $X$ \cite{waldy}. 
For example, the $A$-theory of a smooth manifold $M$ is related to a space of stable $h$-cobordisms which contains information about the diffeomorphism group of $M$, via the stable
parametrized $h$-cobordism theorem \cite{waldy, WBR}.
When $X$ is connected, then $A(X)$ is equivalent to the algebraic $K$-theory of the ring $\SI \Omega X$.
In \cite{HScomod}, Hess--Shipley provide a coalgebraic framework for variants of $A$-theory as algebraic $K$-theories determined by homotopically finite $X_+$-comodules in pointed spaces. 
Similarly, in \cite{perDoldKan}, the second author further proved that a rational variant of $A(X)$ is equivalent to the algebraic $K$-theory of  comodules over the singular chain complex $C_*(X;\mathbb{Q})$ that are perfect differential graded chain complexes.
These algebraic $K$-theories of comodules, arising out of $A$-theory, more closely resemble coalgebraic versions of classical $G$-theories  rather than the coalgebraic $K$-theory $K^c$ discussed above. 

Given a field $k$ and a $k$-algebra $A$, the $G$-theory of $A$, denoted $G^k(A)$, is closely related to
the algebraic $K$-theory of $A$, but based on $A$-modules that are finite-dimensional as $k$-vector spaces.
More explicitly, $G^k(A)$ is the $K$-theory of finite-dimensional representations of $A$. 
In \cite{cotrace}, the second author and Klanderman define an analogous $G$-theory for a $k$-coalgebra $C$, that we denote here $G^c(C)$, determined by $C$-comodules that are finite dimensional as vector spaces.
We show in this paper how this coalgebraic $G$-theory is related to the classical $G^k$-theory recalled above. 

\begin{introtheorem}[\cref{corollary: equivalence G(A)=G(A^*)}]
Let $k$ be a field and $A$ a $k$-algebra. There exists an equivalence of algebraic $K$-theory spectra:
        \[
        \bK(A^\circ)\simeq G^k(A).
        \]
        It is an  equivalence of ring spectra if $A$ is a bialgebra, and an equivalence of commutative ring spectra if $A$ is a cocommutative bialgebra. 
\end{introtheorem}

More generally, we establish a comparison map $G^c(C)\to G^k(C^*)$ for any $k$-coalgebra $C$, that is an equivalence under mild hypotheses on $C$ 
(\cref{prop: KP identification of co-G-theory} and \cref{cor: equivalence of G^c(C) with G(C^*)}). 

For every $k$-algebra $A$, there is a group homomorphism $\chi\colon G^k_0(A)\to A^*$ induced by the character of finite-dimensional representations.
In \cite{cotrace}, Klanderman and the second author also establish a group homomorphism $G^c_0(C)\to \coHH_0(C)$  for any $k$-coalgebra $C$, that we denote by $\chi^c$ here.
We show in \cref{prop: compatibility between characters} and \cref{corollary: equivalence G(A)=G(A^*)} that the comparison maps between $G^c$-theory and $G^k$-theory above are compatible with the homomorphisms $\chi^c$ and $\chi$.
We apply this to the particular case where $A$ is a group ring, $k\Gamma$, for a group $\Gamma$. The finite dual coalgebra of $k\Gamma$ is given by the so-called representative bialgebra $R_k(\Gamma)$ defined as the subvector space of $\Map(\Gamma, k)$ spanned by the representative functions. This allows us to give a coalgebraic interpretation of the classical character homomorphism of a group.

\begin{introtheorem}[\cref{cor: Swan theory as coalgebras}]
     Given any group $\Gamma$, there is an equivalence of algebraic $K$-theory spectra:
    \[
    \bK(R_k(\Gamma))\simeq G^k(k\Gamma).
    \]
    Moreover, the character  $\chi\colon G^k_0(k\Gamma)\rightarrow \Map(\Gamma, k)$ is recovered by the trace $\chi^c\colon \bK_0(R_k(\Gamma))\to\coHH_0(R_k(\Gamma))$.
\end{introtheorem}

\subsection*{Organization}
In \cref{section: intro of K^c and G^c}, we recall the analogues of algebraic $K$-theory and $G$-theory in the coalgebraic setting. In particular, we recall the definitions of $K^c(C)$ and $G^c(C)$ for a coalgebra $C$ over a field $k$.
In \cref{sec:coreflexive}, we recall the definition of coreflexive coalgebras and provide examples.
In \cref{sec: comparison map between K^c and K}, we construct a comparison map between coalgebraic $K$-theory and classical algebraic $K$-theory, $K^c(C)\to K(C^*)$, and show that it is compatible with (co)trace maps. In \cref{sec: comparison between G^c and G} we construct a  comparison map between coalgebraic $G$-theory and classical $G$-theory, $G^c(C)\to G^k(C^*)$, and demonstrate compatibility with character maps. 
In \cref{sec: comparison on Sweedler duals}, we study these comparison maps for the Sweedler dual coalgebra $A^\circ$ of a $k$-algebra $A$, and further determine conditions under which the comparison maps are equivalences.
Lastly, in \cref{sec: app to swan theory}, we apply the results of the previous section to the group algebra $k\Gamma$ of a group $\Gamma$, and provide a coalgebraic framework for the usual character of a group representation.

\subsection*{Acknowledgments:}
The authors are grateful to David Chan for helpful conversations related to this work. The authors also thank the anonymous referee for their helpful comments. The first author was supported by NSF grants DMS-2104233 and DMS-2404932. This work was also partially supported by a grant from the Simons Foundation. The first and second authors would like to thank the Isaac Newton Institute for Mathematical Sciences, Cambridge, for support and hospitality during the programme ``Equivariant homotopy theory in context" where work on this paper was undertaken. This work was supported by EPSRC grant no EP/Z000580/1. The third author was supported by a Fulbright grant for a research stay at Michigan State University.
This research was further supported by NSF grant DMS-2135960.

\section{Coalgebras and their algebraic \texorpdfstring{$K$}{TEXT}-theories}\label{section: intro of K^c and G^c}

In this section we recall some basic definitions for coalgebras, as well as analogues of algebraic $K$-theory and $G$-theory in the coalgebraic setting.  Throughout, let $k$ be a field and $\Vect_k$ denote the category of $k$-vector spaces and $k$-linear homomorphisms.

\begin{definition}
Given a field $k$, recall that a \textit{$k$-coalgebra} is a comonoid in the symmetric monoidal category of $k$-vector spaces with their tensor product $\otimes_k=\otimes$, i.e.~a monoid (or algebra object) in $\Vect_k^\op$.
In other words, a $k$-coalgebra $(C, \Delta, \varepsilon)$ consists of a $k$-vector space $C$, together with $k$-linear homomorphisms $\Delta\colon C\to C\otimes C$ (called the comultiplication) and $\varepsilon\colon C\to k$ (called the counit), that are coassociative and counital. This means $(\id_C\otimes \Delta)\circ \Delta = (\Delta\circ\id_C)\circ \Delta$ and $(\id_C\otimes \varepsilon)\circ \Delta=\id_C=(\varepsilon\otimes \id_C)\circ \Delta$.
We may refer to the triple $(C,\Delta, \varepsilon)$ just as $C$ if the comultiplication and counit are understood.
We shall occasionally use the Sweedler notation and write $\Delta(c)=\sum_i c_{(1)_i}\otimes c_{(2)_i}$ simply as:
\[
\Delta(c) =\sum_{(c)}c_{(1)}\otimes c_{(2)}\in C\otimes C.
\]
A coalgebra $C$ is \textit{cocommutative} if $\tau\circ \Delta=\Delta$, where $\tau\colon C\otimes C\to C\otimes C$ swaps the terms, i.e., $\tau(c\otimes c')=c'\otimes c$. In other words, $C$ is cocommutative if for all $c\in C$:
\[
\sum_{(c)}c_{(1)}\otimes c_{(2)}=\sum_{(c)}c_{(2)}\otimes c_{(1)}.
\]
A \textit{homomorphism of coalgebras} $(C,\Delta_C, \varepsilon_C)\to (D, \Delta_D, \varepsilon_D)$ consists of a $k$-linear homomorphism $f\colon C\to D$ such that $\Delta_D\circ f=(f\otimes f)\circ \Delta_C$ and $\varepsilon_D\circ f=\varepsilon_C$. We denote by $\coAlg_k$ the induced category of $k$-coalgebras with coalgebra homomorphisms. Notice that $\coAlg_k$ is equivalent to $\Alg(\Vect_k^\op)$, the category of algebra objects in $\Vect_k^\op$.
\end{definition}

\begin{definition}
    A \textit{$k$-bialgebra} $H$ is a $k$-vector space $H$ together with a $k$-coalgebra structure $(H, \Delta, \varepsilon)$ and a $k$-algebra structure $(H, \mu, \eta)$ such that $\Delta\colon H\to H\otimes H$ and $\varepsilon\colon H\to k$ are algebra homomorphisms (or equivalently, $\mu\colon H\otimes H\to H$ and $\eta\colon k\to H$ are coalgebra homomorphisms).
    A bialgebra $H$ is \textit{commutative} if it is commutative as an algebra, and $H$ is \textit{cocommutative} if it is cocommutative as a coalgebra. 
    A bialgebra $H$ is a \textit{Hopf algebra} if there exists a $k$-linear function $S\colon H\to H$ (necessarily unique) such that $\mu\circ (\id\otimes S)\circ \Delta = \eta \circ \varepsilon = \mu\circ (S\otimes \id)\circ \Delta$.
\end{definition}

\begin{definition}
Given a $k$-coalgebra $C$, a \textit{right $C$-comodule} $(M, \rho)$ consists of a $k$-vector space $M$ together with a $k$-linear homomorphism $\rho\colon M\to M\otimes C$ that is coassociative and counital: $(\id_M\otimes \Delta)\circ \rho=(\rho\otimes \id_C)\circ \rho$ and $(\id_M\otimes \varepsilon)\circ \rho=\id_M$.
A (right) \textit{$C$-colinear homomorphism} $f\colon (M,\rho)\rightarrow (M', \rho')$ is a $k$-linear homomorphism $f\colon M\rightarrow M'$ such that $\rho'\circ f= (f\otimes \id_C)\circ \rho$.
Let $\coMod_C$ denote the category of right $C$-comodules with colinear homomorphisms.
\textit{Left $C$-comodules} are defined completely analogously. 
If $C$ is cocommutative, then  left and right $C$-comodules  are equivalent and we will simply refer to them as $C$-comodules. We shall occasionally use the Sweedler notation for the coaction and write simply $\sum_{(m)}m_{(0)}\otimes m_{(1)}$ for $\rho(m)=\sum_{i} {m_{(0)}}_i\otimes {m_{(1)}}_i\in M\otimes C$ for any $m\in M$.
\end{definition}

\begin{definition}
A right $C$-comodule $M$ is \textit{finitely cogenerated} if there exists a $C$-colinear monomorphism $M\hookrightarrow C^{\oplus n}\coloneqq k^{\oplus n}\otimes C$. 
\end{definition}

\begin{definition}
A right $C$-comodule $M$ is \textit{injective} if for every $C$-colinear monomorphism $\iota\colon X\hookrightarrow Y$ and any $C$-colinear homomorphism $f\colon X\to M$, there exists a $C$-colinear homomorphism $g\colon  Y\to M$ such that $g\circ \iota=f$.
\end{definition}

In particular, if $M$ is a finitely cogenerated and injective right $C$-comodule, there exists another finitely cogenerated and injective right $C$-comodule $N$ such that $M\oplus N\cong C^{\oplus n}$ as comodules, for some $n\geq 0$. Let $\Inj_\fc(C)$ denote the category of finitely cogenerated and injective right $C$-comodules. 

\begin{proposition}
The category $\Inj_\fc(C)$ is an exact category. 
\end{proposition}

\begin{proof}
    The category $\coMod_C$ is an abelian category as finite limits and colimits in $\coMod_C$ are created under the forgetful functor $\coMod_C\to \Vect_k$.
    Consider a short exact sequence of right $C$-comodules:
    \[
    \begin{tikzcd}
        0 \ar{r} & M \ar{r} & N \ar{r} & P \ar{r} & 0.
    \end{tikzcd}
    \]
    If $M$ and $P$ are finitely cogenerated and injective, then so is $M\oplus P$ and the short exact sequence splits, thus as $N\cong M\oplus P$, we can conclude.
\end{proof}

\begin{definition}[{\cite{cotrace}}]
    Given a $k$-coalgebra $C$, its \textit{coalgebraic $K$-theory} $K^c(C)$ is the algebraic $K$-theory spectrum $K(\Inj_\fc(C))$ of the exact category $\Inj_\fc(C)$ of finitely cogenerated and injective right $C$-comodules.
\end{definition}

The class of finitely cogenerated and injective comodules forms the class of dualizable objects in comodules with respect to a monoidal structure we now make precise.
Recall that given a right $C$-comodule $(M, \rho)$ and a left $C$-comodule $(N, \lambda)$, the relative cotensor product $M\square_C N$ is defined as the equalizer in $\Vect_k$:
\[
\begin{tikzcd}
M\square_C N \ar{r} & M\otimes N \ar[shift left]{r}{\rho\otimes 1} \ar[shift right]{r}[swap]{1\otimes \lambda} & M\otimes C \otimes N.
\end{tikzcd}
\]
In other words, we can view $M\square_C N$ as the relative tensor product in the category of modules over the ``algebra" $C$ in  $\Vect_k^\op$. 
Because the tensor product in $\Vect_k$ preserves limits, the vector space $M\square_C N$ is a $C$-comodule whenever $C$ is cocommutative.
In that case, $\coMod_C$ is a symmetric monoidal category with respect to the cotensor product $\square_C$ with unit $C$, and the dualizable comodules are precisely $\Inj_\fc(C)$ \cite{cotrace}.
More generally, if $C$ is not cocommutative, the bicomodules over $C$ form a (nonsymmetric) monoidal structure with the cotensor product, and we view finitely cogenerated and injective comodules as certain dualizable objects in a bicategory of bicomodules \cite{cotrace}.

Recall that if $R$ is a Noetherian ring, the $G$-theory of $R$, denoted $G(R)$, is the algebraic $K$-theory of the category of finitely generated right $R$-modules (see, for instance, \cite{Quillen, weibel}). These are precisely the compact generators in the category of right $R$-modules. 
The fundamental theorem of comodules shows that, when considering $C$ a coalgebra over a field $k$, every right $C$-comodule is the colimit of its finite dimensional subcomodules \cite{Sweedler}.
In particular, finite dimensional comodules assemble into an abelian category denoted $\coMod_C^\fd$ and they form the compact generators of $\coMod_C$. 
This yields the following analogue of $G$-theory in the coalgebraic setting.

\begin{definition}[{\cite{cotrace}}]
    Given a $k$-coalgebra $C$, the \textit{coalgebraic $G$-theory $G^c(C)$ of $C$} is the algebraic $K$-theory spectrum $K(\coMod_C^\fd)$ of the abelian category $\coMod_C^\fd$ of right $C$-comodules that are finite dimensional as vector spaces.
\end{definition}

As defined  by Swan, the algebraic $G$-theory of a $k$-algebra $A$, denoted $G^k(A)$, is the algebraic $K$-theory spectrum of the abelian category $\Mod_A^\fd$ of right $A$-modules that are finite dimensional as vector spaces \cite{swan}.
We shall see below that $G^c(C)$ is more closely related to $G^k(A)$ than to $G(A)$.
When $A$ is a finite dimensional $k$-algebra, then $G^k(A)$ is simply $G(A)$. A particular example of interest is when $A=k\Gamma$, the group algebra of a group $\Gamma$. In this case, $G^k(k\Gamma)$  is referred to as Swan theory \cite{swan}. See \cref{sec: app to swan theory} for more details. 

\section{(Co)reflexive (co)algebras}\label{sec:coreflexive}
In the current work, an important class of coalgebras will be coreflexive coalgebras. In this section, we briefly recall the definitions of reflexive algebras and coreflexive coalgebras from \cite{Taf72}.

The linear dual functor $(-)^*\colon \Vect_\FF^\op\rightarrow \Vect_\FF$ enjoys the following properties.
\begin{enumerate}
    \item It is a right adjoint, with left adjoint given by itself ${(-)^*}^\op\colon \Vect_\FF\rightarrow \Vect_\FF^\op$.
    It induces an equivalence of symmetric monoidal categories when restricted to finite dimensional vector spaces:
    \[
    (-)^*\colon (\Vect_\FF^\fd)^\op\stackrel{\cong}\longrightarrow \Vect_\FF^\fd.
    \]
    \item It is an exact functor as $\FF$ is injective as a $\FF$-module and thus $\Hom(-,\FF)=(-)^*$ is exact.
    \item It is lax symmetric monoidal, induced by the natural inclusion for any vector spaces $V$ and $W$:
\begin{align*}
    \Phi_{V,W}\colon V^*\otimes W^* & \hookrightarrow (V\otimes W)^*\\
    f \otimes g & \mapsto \left( \begin{array}{rl}
V\otimes W & \rightarrow \FF \\
v\otimes w & \mapsto f(v)g(w)
\end{array}\right).
\end{align*}
If either $V$ or $W$ is finite dimensional, then $\Phi_{V,W}$ is an isomorphism.
\end{enumerate}

A first well-known consequence of the linear dual being lax symmetric monoidal is that it preserves algebra objects, and thus yields a functor $(-)^*\colon \mathrm{Alg}(\Vect_\FF^\op)\rightarrow \Alg(\Vect_\FF).$ This equivalently  defines a functor: 
\[
(-)^*\colon \coAlg_\FF^\op\rightarrow \Alg_\FF.
\]
This is a reformulation of the classical result that the dual $C^*$ of a $\FF$-coalgebra $C$ is a $\FF$-algebra.
The multiplication on $C^*$ is often referred to as the \textit{convolution product}. Explicitly, if we denote the comultiplication $\Delta\colon C\rightarrow C\otimes C$ as $\Delta(c)=\sum_{(c)}c_{(1)}\otimes c_{(2)}$, then the multiplication $f\star g$ of elements in  $C^*$ is defined as:
\[
(f\star g)(c)=\sum_{(c)}f(c_{(1)})g(c_{(2)})
\]
for all $c\in C$.
Because the linear dual is a symmetric monoidal equivalence on finite dimensional vector spaces, there is in fact an (anti-)equivalence between coalgebras and algebras:
\[
(-)^*\colon (\coAlg_\FF^\fd)^\op \stackrel{\cong}\longrightarrow \Alg_\FF^\fd,
\]
that is its own inverse. Indeed, given a finite dimensional algebra $A$, its dual $A^*$ is a coalgebra. But in general, without the finite dimensional condition, there is no comultiplication on $A$ a priori. 
Without finite dimensionality, there is an adjunction
\begin{equation}
    \begin{tikzcd}[column sep=large]
        \Alg_\FF \ar[bend left]{r}{(-)^\circ}\ar[phantom, "\perp" description, xshift=0ex]{r} & \coAlg_\FF^\op \ar[bend left]{l}{(-)^*}
    \end{tikzcd}
\end{equation}
where $A^\circ$ is called the finite dual (or Sweedler dual) of $A$. This was introduced in \cite{Sweedler} and is defined as
\[
A^\circ\coloneqq \left\lbrace f\in A^* \mid \exists I \trianglelefteq A: \dim(A/I)< \infty \text{ and } I\subseteq \ker(f) \right\rbrace.
\]
The comultiplication on $A^\circ$ is given by:
\begin{align*}
    \Delta\colon A^\circ & \longrightarrow A^\circ\otimes A^\circ=(A\otimes A)^\circ\\
    f & \longmapsto \Delta f
\end{align*}
where $\Delta f(a\otimes b)=f(ab)$, for all  $a,b\in A$. 
Its counit is given by 
\begin{align*}
    \varepsilon\colon A^\circ & \longrightarrow k\\
    f & \longmapsto f(1_A).
\end{align*}
There are natural maps $\varepsilon_C\colon C\rightarrow (C^*)^\circ$ and $\eta_A\colon A\rightarrow (A^\circ)^*$, the counits and units of the adjunctions. In other words, there is a natural identification:
\[
\coAlg_k(C, A^\circ)\simeq \Alg_k(A, C^*)
\]
for any $k$-algebra $A$ and any $k$-coalgebra $C$. In general, the map $\varepsilon_C\colon C\rightarrow (C^*)^\circ$ is always injective, however the map $\eta_A\colon A\rightarrow (A^\circ)^*$ is neither injective nor surjective in general. We recall the following definition from \cite{Taf72}. 

\begin{definition}
A $k$-coalgebra $C$ is called \textit{coreflexive} if the natural coalgebra homomorphism $\varepsilon_C\colon C\rightarrow (C^*)^\circ$ is an isomorphism.
A $k$-algebra $A$ is \textit{reflexive} if the natural algebra homomorphism $\eta_A\colon A\rightarrow (A^\circ)^*$ is an isomorphism. 
The $k$-algebra $A$ is \textit{proper} if $\eta_A$ is only assumed to be injective, and $A$ is \textit{weakly-reflexive} if $\eta_A$ is only assumed to be surjective.
\end{definition}

Essentially, asking a $k$-algebra $A$ to be reflexive is requiring $A$ to be obtained from a coalgebra $C$ in the sense that $C^*=A$.
More precisely, by duality, a coalgebra $C$ is coreflexive if and only if $C\cong A^\circ$ for some reflexive algebra $A$; and $A$ is a reflexive algebra if and only if $A\cong C^*$ for some coreflexive coalgebra $C$ \cite[1.8, 1.9]{Rad73}. 
Of course, if $A$ or $C$ are finite dimensional, then they are reflexive or coreflexive respectively. There are also non-finite dimensional examples.

\begin{example}\label{example: k[x] is coreflexive}
We denote by $\dkx$ the $k$-algebra of divided powers.
Recall that it is a countable vector space generated by $X^{[0]}=1, X^{[1]}, X^{[2]}, \ldots$ with multiplication given by:
\(
X^{[i]}X^{[j]}=\binom{i+j}{i}X^{[i+j]}.
\)
It is a $k$-coalgebra with comultiplication $\Delta(X^{[n]})=\sum_{i=0}^n X^{[i]}\otimes X^{[n-i]}$ and counit $\varepsilon(X^{[n]})=0$ for $n\geq 1$, and $\varepsilon(1)=1$.
In fact, it becomes a Hopf algebra with antipode $S(X^{[n]})=(-1)^n X^{[n]}$. 
When $\mathrm{char}(k)=0$, there is an isomorphism of $k$-algebras with the usual polynomial ring $k[x]$:
\begin{align*}
    \dkx & \stackrel{\cong}\longrightarrow k[x]\\
    X^{[n]} & \longmapsto \frac{x^n}{n!}.
\end{align*}
In fact, the polynomial ring $k[x]$ becomes also a Hopf algebra with comultiplication
$\Delta(x^n)=\sum_{i=0}^n\binom{n}{i}x^{i}\otimes x^{n-i}$, counit $\varepsilon(x^n)=0$ if $n\geq 1$ and $\varepsilon(1)=1$, and antipode $S(x^n)=(-1)^nx^n$. If  $\mathrm{char}(k)=0$, then the isomorphism of algebras $\dkx \cong k[x]$ is actually an isomorphism of Hopf algebras with respect to these comultiplications. 
Considering the linear dual, there is an isomorphism of $k$-algebras:
\begin{align*}
    (\dkx)^* & \stackrel{\cong}\longrightarrow k[[y]]\\
    f & \longmapsto \sum_{n\geq 0} f(X^{[n]})y^n
\end{align*}
with the $k$-algebra of formal power series. Similarly, there is the identification of $k$-algebras $(k[x])^*\cong \dky$, the $k$-algebra of formal divided power series, with multiplication $Y^{[i]}Y^{[j]}=\binom{i+j}{i}Y^{[i+j]}$. 
Again, if $\mathrm{char}(k)=0$, then $\dky\cong k[[y]]$ as $k$-algebras.
We can identify $k[[y]]^\circ$ via an isomorphism of coalgebras 
\begin{align*}
    (k[[y]])^\circ & \stackrel{\cong}\longrightarrow \dkx\\
    f & \longmapsto \sum_{i=0}^{n-1} f(y^i)X^{[i]}
\end{align*}
where $n$ is chosen to be the largest natural number such that $(y^n)$ is contained in $\ker(f)$. Similarly $(\dky)^\circ\cong k[x]$ as $k$-coalgebras.
Therefore $k[x]$ and $\dkx$ are coreflexive $k$-coalgebras, and $k[[y]]$ and $\dky$ are reflexive $k$-algebras.
\end{example}

\begin{example}\label{ex: dual of C Noetherian implies C coreflexive}
    More generally, if $C^*$ is a commutative Noetherian ring, then $C$ is a coreflexive coalgebra \cite[2.2.16]{abe}. 
\end{example}

\begin{example}
    As $k$-algebras, the rings $k[x]$ and $\dkx$ are \emph{not} reflexive $k$-algebras.
    The $k$-coalgebra $\Lin(k)$ is the finite dual $k[x]^\circ$ and is referred as the Hopf algebra of recursive linear functions \cite{recursive}. 
\end{example}

\begin{example}[{\cite[\S 3]{Taf72}}]
The dual $C^*$ of a coalgebra is always proper. Finitely generated commutative algebras, free algebras, and free commutative algebras are always proper.
\end{example}

If $A$ is a weakly-reflexive algebra, then $A^\circ$ is a coreflexive coalgebra \cite[6.2]{Taf72}. However, this condition is not necessary, as demonstrated by the following example.

\begin{example}[{\cite[\S 6]{Taf72}}]
    The polynomial ring $A=\overline{\mathbb{Q}}[x]$ is not a weakly-reflexive algebra, but its finite dual $A^\circ=\Lin(\overline{\mathbb{Q}})$ is a coreflexive coalgebra.
\end{example}

\section{A comparison map between coalgebraic \texorpdfstring{$K$}{TEXT}-theory and algebraic \texorpdfstring{$K$}{TEXT}-theory} \label{sec: comparison map between K^c and K}

In this section, given a $k$-coalgebra $C$, we establish a comparison map between the coalgebraic $K$-theory $K^c(C)$ of $C$ introduced in \cref{section: intro of K^c and G^c}, and the classical algebraic $K$-theory $K(C^*)$ of the dual $k$-algebra $C^*$ seen as a ring.

As the linear dual $(-)^*\colon \Vect_k^\op\rightarrow \Vect_k$ is lax symmetric monoidal, it induces a right adjoint functor
\[
(-)^*\colon \coMod_C^\op \longrightarrow \Mod_{C^*}.
\]
Explicitly, given a right $C$-comodule $M$ with coaction $\rho\colon M\rightarrow M\otimes C$, $M^*$ becomes a right $C^*$-module:
\[
\begin{tikzcd}
    M^*\otimes C^* \ar[hook]{r}{\Phi} & (M\otimes C)^* \ar{r}{\rho^*} & M^*.
\end{tikzcd}
\]
In other words, given $f\in M^*$ and $\alpha\in C^*$, the right action $f\cdot \alpha$ is determined as:
\[
(f\cdot \alpha)(m)= \sum_{(m)} f(m_{(0)})\alpha(m_{(1)})
\]
where $\rho(m)=\sum_{(m)}m_{(0)}\otimes m_{(1)}$ for all $m\in M$.

\begin{remark}\label{remark: formal linear dual comodule is module}
    Formally, the functor $(-)^*\colon \coMod_C^\op \longrightarrow \Mod_{C^*}$ is obtained as follows. Let  $\cC$ and $\cD$ be symmetric monoidal categories.
    Let $R\colon \cC\rightarrow \cD$ be a lax symmetric monoidal right adjoint functor.
    Let $A$ be an algebra in $\cC$. Then $R$ lifts to a functor on the module categories $R\colon \Mod_A(\cC)\rightarrow \Mod_{R(A)}(\cD)$. Apply this discussion when $R$ is the linear dual $(-)^*$ and notice that $\coMod_C^\op=\Mod_C(\Vect_k^\op)$.
\end{remark}

If a right $C$-comodule $M$ is finitely cogenerated, then $M^*$ is finitely generated as a right $C^*$-module: a $C$-colinear monomorphism $M\hookrightarrow C^{\oplus n}$ corresponds to a $C^*$-linear epimorphism $(C^*)^{\oplus n}\rightarrow M^*$.
If $M$ is finitely cogenerated and injective, then there exists a right $C$-comodule $N$ such that $M\oplus N\cong C^{\oplus n}$, for some $n\geq 0$.
Thus $M^*\oplus N^*\cong (C^*)^{\oplus n}$, and so $M^*$ is also a projective $C^*$-module.

Let $\Inj_\fc(C)$ be the category of finitely cogenerated and injective right $C$-comodules. Therefore there is a functor:
\[
(-)^*\colon \Inj_\fc(C)^\op\longrightarrow \Proj_\fg(C^*)
\]
which preserves exact sequences, as $[-, k]\colon \Vect_k^\op\rightarrow \Vect_k$ is exact.
If $C$ is finite dimensional, the functor above is an equivalence of exact categories. 

Recall from \cite[\S 2, Equation (7)]{Quillen} that if $\mathcal{A}$ is an exact category, then $K(\mathcal{A})\simeq K(\mathcal{A}^\op)$. 
Therefore we have obtained the following.

\begin{theorem}\label{thm coalg KT to alg KT}
    Let $C$ be a $k$-coalgebra.
    There exists a map of algebraic $K$-theory spectra:
    \[
    \cK(C)\longrightarrow K(C^*)
    \]
    that becomes an equivalence if $C$ is finite dimensional.
\end{theorem}

Recall from \cref{section: intro of K^c and G^c} that if $C$ is a cocommutative $k$-coalgebra, then the category of right $C$-comodules is symmetric monoidal with respect to the cotensor product $\square_C$. As the cotensor product of finitely cogenerated and injective comodules remains a finitely cogenerated and injective comodule \cite[2.20]{cotrace}, $\Inj_\fc(C)$ is also a symmetric monoidal category with respect to the cotensor product. This shows that $K^c(C)$ is a commutative ring spectrum if $C$ is cocommutative (see for instance \cite[2.8]{blum-man}) when viewing $\Inj_\fc(C)$ as a Waldhausen category with weak equivalences given by isomorphisms, and cofibrations by monomorphisms.

If $C$ is cocommutative, then its dual $C^*$ is a commutative $k$-algebra.
Thus the usual algebraic $K$-theory spectrum $K(C^*)$ is also a commutative ring spectrum. We next show that \cref{thm coalg KT to alg KT} can be promoted to a ring homomorphism.

\begin{theorem}\label{thm: comparison map of K^c is ring hm}
If $C$ is a cocommutative $k$-coalgebra, then the natural map of algebraic $K$-theory spectra:
    \[
    \cK(C)\longrightarrow K(C^*)
    \]
from \cref{thm coalg KT to alg KT} is a ring homomorphism of commutative ring spectra.
\end{theorem}

\begin{proof}
    It suffices to show that the functor:
    \[
    (-)^*\colon \coMod_C^\op \longrightarrow \Mod_{C^*}
    \]
    is lax symmetric monoidal (see, for instance, \cite[\S 2]{blum-man}).
    The lax symmetric monoidal structure of the linear dual $\Phi_{M, N}\colon M^*\otimes N^*\hookrightarrow (M\otimes N)^*$ lifts to
    a lax symmetric monoidal structure on $C^*$-modules:
    \begin{align*}
       M^*\tens_{C^*} N^* &\longrightarrow \left(M\cotens_C N\right)^*  \\
       f \otimes g & \longmapsto \left(
       \begin{array}{rcl}
       \displaystyle M\square_C N & \rightarrow & k\\
       m\otimes n & \mapsto & f(m)g(n).
       \end{array}\right)
    \end{align*}
    for all right $C$-comodules $M$ and $N$. 
    
    Formally, the natural map above is obtained as follows. 
    Let $(\cC, \wedge, \mathbbm{1}_\cC)$ and $(\cD, \odot, \mathbbm{1}_\cD)$ be symmetric monoidal categories such that their monoidal products preserve coequalizers in each variable.
    Let $R\colon \cC\rightarrow \cD$ be a lax symmetric monoidal right adjoint functor. Following \cref{remark: formal linear dual comodule is module},
if  $A$ is a commutative algebra in $\cC$, then $R$ lifts to a functor on the module categories $R\colon \Mod_A(\cC)\rightarrow \Mod_{R(A)}(\cD)$.
    Because of our hypothesis, the module categories $\Mod_A(\cC)$ and $\Mod_{R(A)}(\cD)$ are symmetric monoidal with respect to their relative tensor products $\wedge_A$ and $\odot_{R(A)}$ respectively. 
    As $R\colon \cC\rightarrow \cD$ is lax symmetric monoidal, then so is $R\colon \Mod_A(\cC)\rightarrow \Mod_{R(A)}(\cD)$.
    We apply this discussion to the linear dual functor $(-)^*\colon \Vect_k^\op\rightarrow \Vect_k$.
\end{proof}

\begin{example}\label{example: what are comodules over k[x]}
Consider the $k$-coalgebra of divided powers $\dkx$ from \cref{example: k[x] is coreflexive}. 
Then, by coassociativity, a right $\dkx$-comodule $M$ must have coactions of the form:
\begin{align*}
    M & \longrightarrow M\otimes \dkx\\
    m & \longmapsto \sum_{n\geq 0} \phi^n(m)\otimes X^{[n]}
\end{align*}
where $\phi\colon M\rightarrow M$ is a locally nilpotent $k$-linear homomorphism, in the sense that given $m\in M$, there exists $N$ large enough such that $\phi^n(m)=0$ for all $n> N$.
Its linear dual $M^*$ has a $k[[y]]$-action given by:
\[
(f\cdot \sum_{n\geq 0} \lambda_ny^n)(m)=\sum_{n\geq 0} \lambda_nf(\phi^n(m))
\] 
for all $f\in M^*$, $m\in M$ and $\lambda_i\in k$. If $M$ if finitely cogenerated and injective, then $M\cong \dkx^{\oplus n}$ for some $n$, i.e.\ is cofree, see \cite[2.6]{leo}.
As $k[[y]]$ is a principal ideal domain, we know that every finitely generated and projective module is also free.
The natural ring homomorphism:
    \[
    \cK(\dkx)\longrightarrow K(k[[y]]),
    \]
is determined by the assignment of sending a comodule $\dkx^{\oplus n}$ to a module $k[[y]]^{\oplus n}$. 
We shall see in \cref{sec: comparison on Sweedler duals} that this is an equivalence of commutative ring spectra.
\end{example}

In the classical setting, the study of algebraic $K$-theory is facilitated by an understanding of Hochschild homology. For instance, given a $k$-algebra $A$, the Hattori--Stallings trace defines a group homomorphism $\tr\colon K_0(A)\to \HH_0(A)$ between the algebraic $K$-theory of $A$ and the Hochschild homology of $A$, which is the zeroth level of the Dennis trace.

Given a finitely generated projective right $A$-module $M$, there exists a module $N$ such that $M\oplus N\cong A^{\oplus n}$. Denote $e^M\colon A^{\oplus n}\to A^{\oplus n}$ the projection $M\oplus N\rightarrow M$.
Denote by $e^M_{ij}$ the composite:
 \[
 \begin{tikzcd}
     A \ar{r}{\iota_j} & A^{\oplus n} \ar{r}{e^M} & A^{\oplus n} \ar{r}{p_i} & A
 \end{tikzcd}
 \]
  where $\iota_j$ is the $j$-th inclusion and $p_i$ is the $i$-th projection. 
  By our choice of $e^M$, there exists $n_M\geq 1$ such that
\[
\tr(M)=\sum_{i=1}^n [e^M_{ii}(1_A)]=\sum_{i=1}^{n_M} [1_A]
\]
where the bracket $[1_A]$ denotes the class of $1_A\in A$ in the quotient image $A\rightarrow A/[A,A]\cong \HH_0(A)$.

 In the coalgebraic setting, there is a classical analogue of Hochschild homology, called \textit{coHochschild homology} \cite{Doi}. We recall the definition of coHochschild homology at the zeroth level here. 

\begin{definition}\label{definition: coHH_0}
    Given a coalgebra $(C, \Delta, \varepsilon)$, its \textit{zeroth coHochschild homology $\coHH_0(C)$} is the kernel of $\Delta- (\tau\circ \Delta)\colon C\to C\otimes C$ where $\tau\colon C\otimes C\to C\otimes C$ swaps the terms. In other words, we can view $\coHH_0(C)$ as the largest vector space in $C$ generated by cocommuting elements.
\end{definition}

In \cite[1.5]{cotrace}, Klanderman and the second author built an analogue of the Hattori--Stallings trace for a $k$-coalgebra $C$, called the \textit{cotrace}, denoted $\cotr\colon K^c_0(C)\to \coHH_0(C)^*$. 
Here $\coHH_0(C)^*$ denotes the linear dual of $\coHH_0(C)$. Given a finitely cogenerated injective right $C$-comodule $M$, there exists a comodule $N$ such that $M\oplus N\cong C^{\oplus n}$. 
Denote by $e^M\colon C^{\oplus n}\rightarrow C^{\oplus n}$ the $C$-colinear projection $M\oplus N\rightarrow M$.
From \cite[2.14]{cotrace}, by our choice of $e^M$, there exists $n_M\geq 1$ such that 
$\cotr(M)\in (\coHH_0(C))^*$ is determined by:
\begin{equation}\label{eq: cotrace formula}
\cotr(M)(c)=\sum_{i=1}^{n_M} \varepsilon(c)
\end{equation}
for any $c\in \coHH_0(C)$, where $\varepsilon\colon C\to k$ is the counit.

The trace and cotraces have similar definitions and they are indeed compatible in a sense we make precise in the following result.

\begin{theorem}\label{theorem: cotrace compatible with K(C^*)}
    Let $C$ be a $k$-coalgebra.
    There exists a natural $k$-linear epimorphism $\HH_0(C^*)\rightarrow (\coHH_0(C))^*$ which fits in the commutative diagram of abelian groups:
    \[
    \begin{tikzcd}
        \cK_0(C) \ar{rr}{(-)^*} \ar{dd}[swap]{\cotr} & & K_0(C^*)\ar{dd}{\tr}\\
        & C^* \ar{dr} \ar{dl} & \\
        \coHH_0(C)^* & & \HH_0(C^*). \ar{ll}
    \end{tikzcd}
    \]
\end{theorem}

\begin{proof}
Recall that the vector space obtained as the commutator $[C^*, C^*]\subseteq C^*$ is generated by elements $[f,g]=f\star g-g\star f\in C^*$ where $\star$ denotes the induced multiplication on $C^*$.
    Given the evaluation homomorphism $C^*\otimes C\rightarrow k$, notice that $[f,g]\otimes c\in C^*\otimes C$ is in the kernel of the evaluation for all choices $f,g\in C^*$ if and only if $c\in \coHH_0(C)$. Indeed
    \[
    [f,g](c)=(f\star g)(c)-(g\star f)(c)=\sum_{(c)}f(c_{(1)})g(c_{(2)})-\sum_{(c)}g(c_{(1)})f(c_{(2)}).
    \]
    The universal property of the quotient yields a $k$-linear homomorphism $\HH_0(C^*)\otimes \coHH_0(C)\rightarrow k$, and thus by adjointness a $k$-linear homomorphism $\HH_0(C^*)\rightarrow (\coHH_0(C))^*$.
    Essentially, the homomorphism allows us to choose a representative of a class in $\HH_0(C^*)$ and evaluate with an element in $\coHH_0(C)$.
    The inclusion $\coHH_0(C)\subseteq C$ yields a surjective map $C^*\rightarrow \coHH_0(C)^*$ that fits in a commutative diagram:
    \[
    \begin{tikzcd}
        & C^* \ar{dr} \ar{dl} & \\
        \coHH_0(C)^* & & \HH_0(C^*). \ar{ll}
    \end{tikzcd}
    \]
    Thus $\HH_0(C^*)\rightarrow (\coHH_0(C))^*$ is also surjective.

    Let $M$ be a finitely cogenerated and injective right $C$-comodule. Choose a retract of $M\hookrightarrow C^{\oplus n}$ as above to define a $C$-colinear homomorphism $e^M\colon C^{\oplus n}\rightarrow C^{\oplus n}$ whose image is $M$. 
    Our choice also determines a splitting of the surjection $(C^*)^{\oplus n}\rightarrow M^*$, that we denote $e^{M^*}\colon (C^*)^{\oplus n}\rightarrow (C^*)^{\oplus n}$.
    By our discussion above:
    \[
    \tr(M^*)(c)=\left( \sum_{i=1}^{n_{M}} [1_{C^*}]\right)(c)=\sum_{i=1}^{n_M} 1_{C^*}(c)=\sum_{i=1}^{n_M}\varepsilon(c),
    \]
    for all $c\in \coHH_0(C)$. Hence, the diagram in the statement of the theorem commutes. 
\end{proof}

\begin{remark}\label{Rem: higher trace maps}
Recall from \cite{Doi} or \cite[4.3]{cotrace} that $\coHH_n(C)$ is  $\mathrm{coTor}^n_{C\otimes C^\op}(C,C)$, where $C^\op$ denotes the coalgebra $C$ with reverse comultiplication.
Given a $k$-algebra $C$, there are natural $k$-linear homomorphisms $\HH_n(C^*)\to \coHH_n(C)^*$ for all $n\geq 0$ extending the one appearing in \cref{theorem: cotrace compatible with K(C^*)}.
Therefore, combining the comparison map $K^c_n(C)\to K_n(C^*)$ with the Dennis trace $K_n(C^*)\to \HH_n(C^*)$ leads to a trace map $K^c_n(C)\to \coHH_n(C)^*$ for all $n\geq 0$.
In upcoming work of Agarwal, Mehrle, and the second author, the trace map $K^c_n(C)\to \coHH_n(C)^*$  is studied further. 
\end{remark}

\section{A comparison map between coalgebraic \texorpdfstring{$G$}{TEXT}-theory and algebraic \texorpdfstring{$G$}{TEXT}-theory}\label{sec: comparison between G^c and G}

In this section, given a $k$-coalgebra $C$, we establish a comparison map between the coalgebraic $K$-theory $G^c(C)$ of $C$ introduced in \cref{section: intro of K^c and G^c}, and the algebraic $K$-theory  $G^k(C^*)$ of finite dimensional representations of the dual $k$-algebra $C^*$.

Given a right $C$-comodule $M$, with coaction $\rho\colon M\rightarrow M\otimes C$, there is a left $C^*$-module action on $M$, denoted $\lambda\colon C^*\otimes M\rightarrow M$,  via the composition:
\[
\begin{tikzcd}
    C^*\otimes M \ar{r}{1\otimes \rho} & C^*\otimes M\otimes C\cong C^*\otimes C\otimes M \ar{r}{\mathrm{ev}\otimes 1} & k\otimes M\cong M. 
\end{tikzcd}
\]
In other words, if we write $\rho(m)=\sum_{(m)} m_{(0)}\otimes m_{(1)}$, then the left action on $M$ is defined for all $f\in C^*$:
\begin{equation}\label{equation: actions induced by coaction}
f\cdot m = \sum_{(m)}f(m_{(1)})m_{(0)}.
\end{equation}
One can check that $M\rightarrow N$ is a homomorphism of right $C$-comodules if and only if it is a homomorphism of left $C^*$-modules. 
This defines a fully faithful exact functor 
\[
\begin{tikzcd}
    \coMod_C \ar[hook]{r} & {}_{C^*}\Mod.
\end{tikzcd}
\]
The functor is an equivalence of categories if and only if $C$ is finite dimensional \cite[4.7]{coringscomodules}.
Since this functor restricts to an exact functor $\coMod_C^\fd\rightarrow {}_{C^*}\Mod^\fd$, the definitions of $G^c(C)$ and $G^k(C^*)$  imply the following equivalence. 

\begin{proposition}\label{prop: KP identification of co-G-theory} 
    Let $C$ be a  $k$-coalgebra.
    There is a natural map of algebraic $K$-theory spectra: 
    \[
    \bK(C)\longrightarrow G^k(C^*),
    \]
    which becomes an equivalence if $C$ is finite dimensional.
\end{proposition}

Although $\coMod_C\hookrightarrow {}_{C^*}\Mod$ is not an equivalence if $C$ is not finite dimensional, one can study when the restriction $\coMod_C^\fd\hookrightarrow{}_{C^*}\Mod^\fd$ is an equivalence for non-finite dimensional coalgebras $C$. 
In general, not every left $C^*$-module needs to be a right $C$-comodule, because not every action is as in \cref{equation: actions induced by coaction}.
A left $C^*$-module $M$ with action $\lambda\colon C^*\otimes M\rightarrow M$ defines an adjoint map $\tilde{\lambda}\colon M\rightarrow [C^*, M]\coloneqq \Hom_k(C^*, M)$, where $\tilde{\lambda}(m)(f)=\lambda(f\otimes m)$ for all $f\in C^*$ and $m\in M$.
When $C$ or $M$ is finite dimensional, the natural inclusion:
\begin{align*}
    M\otimes C & \hookrightarrow   [C^*, M]\\
    m\otimes c & \mapsto   \left( \begin{array}{ccl}
C^* & \rightarrow & M\\
f & \mapsto & f(c)m
\end{array}\right).
\end{align*}
is an isomorphism. \cref{equation: actions induced by coaction} is equivalent to the commutative diagram:
\begin{equation}\label{equation: rational modules explanation}
 \begin{tikzcd}
    & M \ar[hook]{dr}{\tilde{\lambda}} \ar[hook']{dl}[swap]{\rho}&\\
    M\otimes C \ar[hook]{rr}{} & & {[C^*, M]}.
\end{tikzcd}   
\end{equation}

\begin{definition}\label{def: rational modules}
    The essential image of $\coMod_C \hookrightarrow {}_{C^*}\Mod$ is often referred to as the \textit{rational $C^*$-modules}.
\end{definition}

Rational modules are related to comodules over coreflexive coalgebras, introduced in \cref{sec:coreflexive}.

\begin{proposition}[{\cite[3.1.3]{abe}}] \label{prop:coreflexive}
Given a $k$-coalgebra $C$, 
all finite dimensional left $C^*$-modules are rational if and only if $C$ is coreflexive.
\end{proposition}

Thus the induced exact  functor $\coMod_C^\fd\hookrightarrow {}_{C^*}\Mod^\fd$ is an equivalence of abelian categories, if and only if $C$ is coreflexive.
Combined with our analysis above, this gives the following identification of $\bK(C)$.

\begin{corollary}\label{cor: equivalence of G^c(C) with G(C^*)}
    If $C$ is a coreflexive $k$-coalgebra, then there is an equivalence of algebraic $K$-theory spectra:
    \[
    \bK(C)\simeq G^k(C^*).
    \]
\end{corollary}

\begin{example}
As $\dkx$ is a coreflexive $k$-coalgebra (\cref{example: k[x] is coreflexive}), there is an equivalence of algebraic $K$-theory spectra:
    \[
    \bK(\dkx)\simeq G^k(k[[y]]).
    \]
    In more detail, if $M$ is a $k[[y]]$-module, then there is a module map $\phi\colon M\rightarrow M$ specified by $\phi(m)=y\cdot m$.
    It is straightforward to check that $M$ is a rational module if and only if $\phi$ is locally nilpotent as in \cref{example: what are comodules over k[x]}.
    This allows us to define a $\dkx$-comodule structure on $M$. 
    Of course, if $M$ is finite dimensional, then $\phi$ must be nilpotent.
\end{example}

Unlike in $\Inj_\fc(C)$, even if $C$ is cocommutative, the cotensor product does not define a monoidal structure on $\coMod_C^\fd$ in general, since the unit $C$ itself might not be finite dimensional as a vector space. Instead, we take the following approach. 
Recall that if $H$ is a bialgebra, then the category of right $H$-comodules $\coMod_H$ is also monoidal with respect to the tensor product $\otimes$ in $\Vect$ instead of the cotensor product. Indeed, given two $H$-comodules $(M, \rho_M)$ and $(N, \rho_N)$, their tensor product $M\otimes N$ becomes an $H$-comodule with coaction: 
\[
\begin{tikzcd}
    M\otimes N \ar{r}{\rho_M\otimes \rho_N} & M\otimes H \otimes N \otimes H\cong M\otimes N\otimes H\otimes H\ar{r}{M\otimes N\otimes \mu} & M\otimes N \otimes H
\end{tikzcd}
\]
where $\mu\colon H\otimes H\rightarrow H$ is the multiplication on $H$. The monoidal unit becomes $k$ itself, seen as an $H$-comodule via the unit $k\rightarrow H$.
If $H$ is commutative, then $\coMod_H$ is a symmetric monoidal category.
An analogous argument works for the modules instead of comodules.
The functor $\coMod_H\hookrightarrow {}_{H^*}\Mod$ is strong monoidal with respect to these monoidal structures whenever $H^*$ is also bialgebra. The following result therefore holds.

\begin{corollary}\label{cor: comparison G^c(C) and G^k(C^*) is monoidal}
    Let $H$ be a $k$-bialgebra.
    Then the natural map $\bK(H)\rightarrow G^k(H^*)$ of \cref{prop: KP identification of co-G-theory} is a map of ring spectra if $H^*$ is also a $k$-bialgebra.
    If moreover $H$ is commutative, then the ring spectra are commutative.
\end{corollary}

In the remainder of this section we show that the comparison map between $G^c(C)$ and $G^k(C^*)$ is compatible with the character homomorphisms on representations.
Let $A$ be a $k$-algebra. 
There is a group homomorphism $\chi\colon G^k_0(A)\rightarrow A^*$ given by the characters, which we now recall.
The data of a finite dimensional right $A$-module is encoded by a $k$-algebra homomorphism $\rho\colon A\rightarrow \mathscr{M}_n(k)$ for some $n$, where $\mathscr{M}_n(k)$ denotes  the $k$-algebra of $n\times n$-matrices with coefficients in $k$. Then $\chi$ is defined as $\chi(\rho)=\tr\circ \rho$, where $\tr\colon \mathscr{M}_n(k)\rightarrow k$ is the usual trace. 
If we choose a basis $(e_1, \ldots, e_n)$ of a finite dimensional right $A$-module $V$ with action $\rho\colon A\rightarrow \mathscr{M}_n(k)$, then for all $a\in A$:
\[
\chi(V)(a)=\sum_{i=1}^n e_i^*(\rho(a)(e_i)).
\]

We show there is a coalgebraic interpretation of this character map.
If $C$ is a $k$-coalgebra, given a finite dimensional right $C$-comodule, the $C$-colinear trace was introduced in \cite[\S 2.3]{cotrace}, and is a group homomorphism $\chi^c\colon \bK_0(C)\rightarrow C$ that factors through $\coHH_0(C)\hookrightarrow C$ and is defined as follows.
Given a finite dimensional right $C$-comodule $V$ with coaction $\rho\colon V\rightarrow V\otimes C$, denoted $v\mapsto \sum_{(v)} v_{(0)}\otimes v_{(1)}$, the character $\chi^c(V)$ is defined as the image of $1$ in the composition:
\[
\begin{tikzcd}
    k \ar{r} & V\otimes V^* \ar{r}{\rho\otimes V^*} & V\otimes C\otimes V^*\cong C\otimes V^*\otimes V \ar{r} & C
\end{tikzcd}
\]
where the first map is the coevaluation on $V$ and the last map is induced by the evaluation on $V$.
If we choose a basis $(e_1, \ldots, e_n)$ of $V$, then there is the formula:
\begin{equation}\label{eq: colinear character}
\chi^c(V)=\sum_{i=1}^n \sum_{(e_i)} e_i^*({e_i}_{(0)}){e_i}_{(1)}.
\end{equation}

\begin{proposition}\label{prop: compatibility between characters}
    Let $C$ be a $k$-coalgebra. 
    Then the character map $\chi\colon G_0^k(C^*)\rightarrow C^{**}$  is compatible with the $C$-colinear trace in the sense that the following diagram commutes:
    \[
    \begin{tikzcd}
        \bK_0(C)\ar{r}{\chi^c} \ar{d} & \coHH_0(C) \ar[hook]{r} & C \ar[hook]{d}{\varepsilon_C}\\
        G^k_0(C^*) \ar{rr}{\chi} & & C^{**}.
    \end{tikzcd}
    \]
\end{proposition}

\begin{proof}
    Let $V$ be a finite dimensional right $C$-comodule, with coaction denoted $v\mapsto \sum_{(v)} v_{(0)}\otimes v_{(1)}$ for all $v\in V$.
    Recall it defines a $C^*$-module via $\rho\colon C^*\rightarrow \mathrm{End}_k(V)$ where 
    \[
    \rho(f)(v)=\sum_{(v)}f(v_{(1)})v_{(0)},
    \]
    for all $v\in V$ and $f\in C^*$.
    Choosing a basis $(e_1, \ldots, e_n)$ of $V$, the character of $\rho$ is then given by:
    \begin{align*}
        \chi(\rho)(f) & = \sum_{i=1}^n e_i^* \left( \rho(f)(e_i)\right)\\
        & = \sum_{i=1}^n e_i^* \left( \sum_{(e_i)} f({e_i}_{(1)}) {e_i}_{(0)}\right)\\
        & = \sum_{i=1}^n \sum_{(e_i)} f({e_i}_{(1)})e_i^*({e_i}_{(0)}),
    \end{align*}
    for all $f\in C^*$.
    On the other hand, recall that for the chosen basis of $V$:
    \[
\chi^c(V)=\sum_{i=1}^n \sum_{(e_i)} e_i^*({e_i}_{(0)}){e_i}_{(1)}.
    \]
    Thus applying the natural embedding $\varepsilon_C\colon C\hookrightarrow C^{**}$ yields:
    \begin{align*}
    \varepsilon_C(\chi^c(V))(f) & = f(\chi^c(V))\\
    & = f\left( \sum_{i=1}^n \sum_{(e_i)} e_i^*({e_i}_{(0)}){e_i}_{(1)}\right)\\
    & = \sum_{i=1}^n \sum_{(e_i)}e_i^*({e_i}_{(0)}) f({e_i}_{(1)}),
     \end{align*}
     for all $f\in C^*$. 
Therefore, the diagram commutes.
\end{proof}

In particular, if the coalgebra $C$ is coreflexive, then $\chi(\rho)\in C^{*\circ}$ for all finite dimensional left $C^*$-modules $\rho\colon C^*\rightarrow \mathrm{End}_k(V)$, since the natural inclusion $\varepsilon_C\colon C\rightarrow C^{**}$ factors through the subvector space $C^{*\circ}$, implying the following result.

\begin{corollary}
    Let $C$ be a coreflexive $k$-coalgebra.
    Then the colinear trace $\chi^c\colon \bK_0(C)\rightarrow \coHH_0(C)$ is a refinement of the character in the sense that the following diagram commutes:
    \[
    \begin{tikzcd}
        \bK_0(C)\ar{r}{\chi^c} \ar{d}[swap]{\cong} & \coHH_0(C) \ar[hook]{r} & C \ar[hook]{d}{\varepsilon_C}[swap]{\cong}\\
        G^k_0(C^*) \ar{rr}{\chi} & & C^{*\circ}.
    \end{tikzcd}
    \]
\end{corollary}

\section{The comparison maps on the Sweedler dual of an algebra}\label{sec: comparison on Sweedler duals}
In \cref{sec: comparison map between K^c and K}, given a $k$-coalgebra $C$, we introduced a comparison map $K^c(C)\to K(C^*)$, and we also obtained a comparison map $G^c(C)\to G^k(C^*)$ in \cref{sec: comparison between G^c and G}.
In this section, we apply and strengthen our previous results in the case where the coalgebra $C$ is obtained as the Sweedler dual $A^\circ$ of an algebra $A$, see \cref{sec:coreflexive}. In this case, we obtain useful reinterpretations of our results in the algebra world.

\begin{definition}
    Let $A$ be a $k$-algebra.
    A left $A$-module $M$ is said to be \textit{locally finite} if for all elements $m$ in $M$, the $k$-vector space $Am$ is finite dimensional.
    In other words, if we denote by $\lambda\colon A\otimes M\rightarrow M$ the left action, the induced linear map $\lambda_m\coloneqq \lambda(-\otimes m)\colon A\rightarrow M$ has finite dimensional rank for all $m\in M$.
\end{definition}

Locally finite $A$-modules form a full subcategory ${}_A\Mod^\lf$ in ${}_A\Mod$. 
The next result can be found in \cite[p. 126]{abe} or \cite[1.7.1]{newtak}.

\begin{proposition}\label{prop: locally finite modules as comodules}
The category of right $A^\circ$-comodules $\coMod_{A^\circ}$ is equivalent to the category ${}_{A}\Mod^\lf$ of locally finite left $A$-modules.
\end{proposition}

We now explain this correspondence.
By extension and restriction  of scalars along the algebra homomorphism $\eta_A\colon A\rightarrow A^{\circ *}$, there is an adjunction:
\[
    \begin{tikzcd}[column sep=large]
        {}_A\Mod \ar[bend left]{r}{A^{\circ *}\otimes_A -}\ar[phantom, "\perp" description, xshift=0.5ex]{r} & {}_{A^{\circ *}}\Mod.\ar[bend left]{l}{}
    \end{tikzcd}
\]
The restriction of scalars is exact. 
As we can see $\coMod_{A^\circ}$ as a full subcategory of rational left $A^{\circ *}$-modules (\cref{def: rational modules}), by combining with restriction of scalars, there is an exact functor
\[
\coMod_{A^\circ} \longrightarrow  {}_A\Mod,
\]
that remains the identity on objects. 
Explicitly, given a right $A^\circ$-comodule $M$ with coaction $M\rightarrow M\otimes A^\circ$ given by $m\mapsto \sum_{(m)} m_{(0)}\otimes m_{(1)}$, the left action $\lambda\colon A\otimes M\rightarrow M$ is defined by:
\[
\lambda(a\otimes m)\coloneqq a\cdot m \coloneqq \sum_{(m)} m_{(1)}(a)m_{(0)}
\]
for all $a\in A$ and $m\in M$.
The image of $\lambda_m\colon A\rightarrow M$ is generated by the elements $m_{(0)}$ which are finite for a fixed choice $m$.
More precisely, given $m$, there are elements ${m_{(0)}}_1, \ldots, {m_{(0)}}_n$ in $M$ for some $n\geq 1$, that generate a subvector space of $M$ that contains $\lambda_m(a)$ for all $a\in A$. Thus $M$ is locally finite as an $A$-module.

Conversely, given a locally finite left $A$-module $M$, with left action $\lambda\colon A\otimes M\rightarrow M$, choosing a basis $\{e_i\}_{i\in \mathcal{I}}$ of $M$, we denote by $\lambda_{ij}(a)\in k$ the scalars obtained as:
\[
\lambda(a\otimes e_i)=\sum_j \lambda_{ij}(a)e_j
\]
for some $a\in A$. Since $M$ is locally finite, the sum over $j$ above is necessarily finite.
Define linear maps $\lambda_{ij}\colon A\rightarrow k$ as $\lambda_{ij}=e_j^*\circ \lambda_{e_i}$.
The kernel of $\lambda_{e_i}\colon A\rightarrow M$ is precisely the torsion of $e_i$, it is an ideal of $A$ and it is cofinite by assumption.
Thus $\lambda_{ij}\in A^\circ$.
We define the right $A^\circ$-coaction on $M$ as:
\begin{align*}
    M & \longrightarrow M\otimes A^\circ\\
    e_i& \longmapsto \sum_{j}e_j\otimes \lambda_{ij}.
\end{align*}
Of course, every finite dimensional left $A$-module is automatically locally finite, therefore we obtain as a consequence the following result.

\begin{corollary}\label{corollary: equivalence G(A)=G(A^*)}
For any $k$-algebra $A$, the functor $\coMod_{A^\circ} \longrightarrow  {}_A\Mod^\lf$ induces an equivalence of abelian categories:
\[
\coMod_{A^\circ}^\fd \cong {}_A\Mod^\fd.
\]
Consequently, there is an equivalence of algebraic $K$-theory spectra: 
\[
\bK(A^\circ)\simeq G^k(A),
\]
that is compatible with the character homomorphism, in the sense that the following diagram commutes: 
\[
\begin{tikzcd}
\bK_0(A^\circ) \ar{r}{\chi^c} \ar{d}[swap]{\cong} & \coHH_0(A^\circ) \ar[hook]{r} & A^\circ \ar[hook]{d}\\
G^k_0(A) \ar{rr}{\chi} & &  A^*
\end{tikzcd}
\]
Moreover, if $A$ is a bialgebra, then $\bK(A^\circ)\simeq G^k(A)$ is an equivalence of ring spectra, and if $A$ is in addition cocommutative, then it is an equivalence of commutative ring spectra.
\end{corollary}

\begin{proof}
The ring structure in the presence of a bialgebra structure follows in a similar way as the discussion above \cref{cor: comparison G^c(C) and G^k(C^*) is monoidal}. 
    We are only left to prove the compatibility of characters.
    Let $V$ be a finite dimensional right $A^\circ$-comodule.
    Let $(e_1, \ldots, e_n)$ be a basis of $V$. 
    Using the Sweedler notation on its coaction, recall that the action of an element $a\in A$ on $v\in V$ is given by:
    \[
    a\cdot v=\sum_{(v)}v_{(1)}(a)v_{(0)}.
    \]
    Therefore the character $\chi(V)\in A^*$ is given by:
    \begin{align*}
        \chi(V)(a) &= \sum_{i=1}^n e_i^*(a\cdot e_i)\\
        & = \sum_{i=1}^n\sum_{(e_i)} e_i^*({e_i}_{(1)}(a){e_i}_{(0)})\\
        & = \sum_{i=1}^n\sum_{(e_i)} {e_i}_{(1)}(a) e_i^*({e_i}_{(0)})
    \end{align*}
    for all $a\in A$. 
    On the other hand, the colinear character $\chi^c(V)$ is given by \autoref{eq: colinear character}:
    \[
    \chi^c(V)(a)=\sum_{i=1}^n\sum_{(e_i)} e_i^*({e_i}_{(0)}) {e_i}_{(1)}(a) 
    \]
    for all $a\in A$, showing that the diagram commutes.
\end{proof}

\begin{example}
    Given a $k[[y]]$-module $M$, we see $M$ is locally finite if the induced map $\phi\colon M\rightarrow M$ defined by $\phi(m)=y\cdot m$ is locally nilpotent as in \cref{example: what are comodules over k[x]}. 
    Therefore, locally finite $k[[y]]$-modules correspond precisely to $\dkx$-comodules.
    Because $k[[y]]$ is reflexive, there is no distinction between locally finite $k[[y]]$-modules and rational $k[[y]]$-modules. 
\end{example}

If $A$ is a $k$-algebra, by \cref{thm coalg KT to alg KT}, the functor $(-)^*\colon (\coMod_{A^\circ})^\op\rightarrow \Mod_{(A^\circ)^*}$ induces a map of algebraic $K$-theory spectra
\[
\cK(A^\circ)\longrightarrow K((A^\circ)^*).
\]
As we argued in \cref{theorem: cotrace compatible with K(C^*)}, the evaluation $A^\circ \otimes A\to k$ sends an element $f\otimes [a,b]$ to zero for all $a,b\in A$ if and only if $f\in \coHH_0(A^\circ)$. 
Indeed, if $f\in \coHH_0(A^\circ)$ then $f(ab)=\Delta f(a\otimes b)=\Delta f(b\otimes a)=f(ba)$, for all $a,b\in A$.
Thus, there is a natural $k$-linear homomorphism $\coHH_0(A^\circ)\otimes \HH_0(A)\to k$ which defines a natural $k$-linear homomorphism:
\[
\HH_0(A) \longrightarrow (\coHH_0(A^\circ))^*.
\]

\begin{theorem}\label{cor: linear dual on coalgebraic K-theory of finite dual}
    For $A$ a reflexive $k$-algebra, there exists a natural map of algebraic $K$-theory spectra:
    \[
    (-)^*\colon \cK(A^\circ)\longrightarrow K(A),
    \]
which fits in the commutative diagram of abelian groups:
    \[
    \begin{tikzcd}
        \cK_0(A^\circ) \ar{rr}{(-)^*} \ar{dd}[swap]{\cotr} & & K_0(A)\ar{dd}{\tr}\\
        & A \ar{dr} \ar{dl} & \\
        \coHH_0(A^\circ)^* & & \HH_0(A). \ar{ll}
    \end{tikzcd}
    \]
    If $A$ is commutative, then $(-)^*\colon \cK(A^\circ)\to K(A)$ is a ring homomorphism of commutative ring spectra.
\end{theorem}

\begin{proof}
    If $A$ is reflexive, then the exact category $\Mod_{(A^\circ)^*}$ is necessarily equivalent to the exact category $\Mod_A$. 
    We can then apply \cref{thm coalg KT to alg KT} and \cref{theorem: cotrace compatible with K(C^*)}.
    The last statement follows from \cref{thm: comparison map of K^c is ring hm}.
\end{proof}

\begin{remark}
    If $A$ is not a reflexive $k$-algebra, then the restriction of scalars $\Mod_{(A^\circ)^*}\rightarrow \Mod_A$ induced by the algebra homomorphism $\eta_A\colon A\rightarrow (A^\circ)^*$ does not in general restrict to finitely generated and projective modules.
    For instance, if $A=k[x]$, then $A^\circ=\Lin(k)$ which has uncountable basis, so $\Lin(k)^*$ is not finitely generated as a $k[x]$-module.
\end{remark}

It can be useful to reinterpret the functor $(-)^*\colon (\coMod_{A^\circ})^\op\rightarrow \Mod_{(A^\circ)^*}$ using the identification of \cref{prop: locally finite modules as comodules}. 
Recall that the linear dual of a left $A$-module $M$ with action $\lambda\colon A\otimes M\rightarrow M$ defines a right $A$-module structure on $M^*$ with action:
\[
\begin{tikzcd}
M^*\otimes A \ar{r}{\lambda^*\otimes A} & (A\otimes M)^*\otimes A \ar[hook]{r} & A^*\otimes M^* \otimes A\cong M^*\otimes A^*\otimes A \ar{r} & M^*
\end{tikzcd}
\]
where the last arrow is induced by evaluation on $A$.
This assignment defines a functor $(-)^*\colon {}_A\Mod^\op \rightarrow \Mod_A$ that fits in the commutative diagram:
\[
\begin{tikzcd}
    \left(\coMod_{A^\circ}\right)^\op \ar{r}{(-)^*} \ar{d}[swap]{\simeq} & \Mod_{A^{\circ *}} \ar{d}\\
    \left( {}_A\Mod^\lf\right)^\op \ar{r}{(-)^*} &  \Mod_A.
\end{tikzcd}
\]
The unlabeled vertical arrow is the restriction of scalars along $\eta_A\colon A\rightarrow A^{\circ *}$.
In general, the linear dual of a locally finite module need not be locally finite. However, there is an adjunction:
\[
\begin{tikzcd}
  \Mod_A  \ar[r, bend left, "(-)^\circ" {xshift=3ex}]
  \ar[phantom, "\perp" description, xshift=1ex]{r} & [1em]\left({}_A\Mod^\lf\right)^\op \ar[phantom, "\simeq", description]{r} \ar[l, bend left, "(-)^*" {xshift=0.5ex}] & [-2 em]\coMod_{A^\circ}^\op
\end{tikzcd}
\]
where for a right $A$-module $M$: \[M^\circ\coloneqq \{f\in M^* \mid \exists N \trianglelefteq M: \dim(M/N)< \infty \text{ and } N\subseteq \ker(f)\}.\]
We can identity $M^\circ$ as the largest locally finite (left) $A$-submodule of $M^*$  \cite[1.7.2]{newtak}.
In general, the left adjoint is only right exact, but it is exact if we restrict to finitely generated right $A$-modules for a commutative Noetherian $k$-algebra $A$ \cite[1.7.4]{newtak}.

Moreover, if $A$ is a commutative $k$-algebra, then $\HH_0(A)=A$, and $A^\circ$ is a cocommutative coalgebra, hence $\coHH_0(A^\circ)=A^\circ$.
This leads to the following result.

\begin{theorem}\label{theorem: finite dual on algebraic K theory of ring}
    Let $A$ be a commutative Noetherian $k$-algebra.
    There exists a natural ring homomorphism of algebraic $K$-theory spectra:
    \[
    (-)^\circ\colon K(A)\longrightarrow \cK(A^\circ)
    \]
    that is compatible with the Hattori--Stallings traces in the sense that the following diagram  of abelian groups commutes:
    \[
    \begin{tikzcd}
    K_0(A) \ar{r}{(-)^\circ} \ar{d}[swap]{\tr}& \cK_0(A^\circ) \ar{d}{\cotr}\\
    A \ar{r}{\eta_A} & A^{\circ *}.
    \end{tikzcd}
    \]
\end{theorem}

\begin{proof}
As $(-)^\circ$ is a left adjoint of an exact functor, it sends a projective $A$-module $M$ to an injective $A$-module $M^\circ$ that is locally finite. As $A$ is Noetherian, it implies that $M^\circ$ is also injective as an $A^\circ$-comodule \cite[1.7.7]{newtak}.
When $M$ is finitely generated as an $A$-module, there is an isomorphism of $A$-modules:
\[
M^\circ \cong  \underline{\Mod}_A(M, A^\circ)
\]
where $\underline{\Mod}_A$ indicates the $A$-module of $A$-linear homomorphisms (i.e.~the internal hom in $\Mod_A$), see \cite[1.7.5]{newtak}.
Since $A^\circ$ is injective as an $A$-module, every surjective $A$-linear homomorphism $A^{\oplus n}\rightarrow M$ gives rise to an injective $A$-linear (and therefore injective $A^\circ$-colinear) morphism :
\[
M^\circ \rightarrow \underline{\Mod}_A(A^{\oplus n}, A^\circ)\cong \underline{\Mod}_A(A, A^\circ)^{\oplus n}\cong (A^\circ)^{\oplus n}.
\]
Thus the functor $(-)^\circ$ sends  finitely generated $A$-modules  to  finitely cogenerated $A^\circ$-comodules,
i.e., it induces an exact functor 
$(-)^\circ\colon \Proj_\fg(A)\rightarrow \Inj_\fc(A^\circ)^\op$.
As we noted for \cref{thm coalg KT to alg KT}, the algebraic $K$-theories of an exact category and its opposite are equivalent.

Since the functor $(-)^\circ$ is a left adjoint, given $A$-modules $M$ and $N$, the relative tensor product $M\otimes_A N$ induces an equalizer of (locally finite) $A$-modules:
\[
\begin{tikzcd}
    (M\otimes A \otimes N)^\circ \ar[shift left, leftarrow]{r} \ar[shift right, leftarrow]{r}  &  (M\otimes N)^\circ \ar[leftarrow]{r} & \displaystyle\left(M\tens_A N\right)^\circ.
\end{tikzcd}
\]
But as $(M\otimes N)^\circ\cong M^\circ\otimes N^\circ$ (see for instance \cite[6.0.1]{Sweedler}), we conclude that the equalizer must be the relative cotensor product of $M^\circ$ with $N^\circ$.  In other words, there is a natural isomorphism of $A^\circ$-comodules:
\[
\displaystyle\left(M\tens_A N\right)^\circ\cong  M^\circ \cotens_{A^\circ} N^\circ.
\]
Therefore $(-)^\circ\colon \Mod_A\rightarrow \coMod_{A^\circ}^\op$ is strong symmetric monoidal.

For the Hattori--Stallings traces, we can argue similarly as in the proof of \cref{theorem: cotrace compatible with K(C^*)}.
Given a finitely generated projective $A$-module $M$, we can choose a splitting of $A^{\oplus n}\to M$ and this induces a choice of a retract of $M^\circ\to (A^\circ)^{\oplus n}$. From this choice of retract, there
exists an integer $n_M$ so that on the one hand: 
\[
\tr(M)=\sum_{i=1}^{n_M}1_A,
\]
and on the other hand:
\[
\cotr(M^\circ)(f)=\sum_{i=1}^{n_M} \varepsilon(f) = \sum_{i=1}^{n_M} f(1_A),
\]
for all $f\in A^\circ$ where $\varepsilon\colon A^\circ\to k$ is the counit. 
As $\eta_A(a)(f)=f(a)$ for all $a\in A$ and $f\in A^\circ$, the diagram commutes, as desired.  
\end{proof}

\begin{corollary}\label{corollary: inverse map K(C^*) to K(C)}
    Let $C$ be a $k$-coalgebra. 
    Suppose its dual $C^*$ is a commutative Noetherian ring.
    Then there exists a natural ring homomorphism of algebraic $K$-theory spectra:
    \[
    (-)^\circ\colon K(C^*)\rightarrow \cK(C),
    \]
    that is compatible with the Hattori--Stallings traces in the sense that the following diagram  of abelian groups commutes:
    \[
    \begin{tikzcd}
    K_0(C^*) \ar{rr}{(-)^\circ} \ar{dr}[swap]{\tr}& & \cK_0(C) \ar{dl}{\cotr}\\
   & C^*.
    \end{tikzcd}
    \]
\end{corollary}

\begin{proof}
By \cref{theorem: finite dual on algebraic K theory of ring}, there is a natural ring homomorphism $K(C^*)\rightarrow \cK(C^{*\circ})$.
Since $C^*$ is commutative Noetherian, $C$ is coreflexive (\cref{ex: dual of C Noetherian implies C coreflexive}).
Therefore $\varepsilon_C\colon C\rightarrow C^{*\circ}$ is an isomorphism of coalgebras, and the result follows.
\end{proof}

Combining \cref{cor: linear dual on coalgebraic K-theory of finite dual} and \cref{theorem: finite dual on algebraic K theory of ring}, there are ring homomorphisms
$(-)^*\colon \cK(A^\circ)\rightarrow K(A)$ and $(-)^\circ\colon K(A)\rightarrow \cK(A^\circ)$ for a nice enough $k$-algebra $A$. Our next result shows these are inverses of each other.

\begin{theorem}\label{theorem: equivalence K(A)=K(A^*)}
    Let $A$ be a commutative Noetherian reflexive $k$-algebra.
    Then the natural maps of \cref{cor: linear dual on coalgebraic K-theory of finite dual} and \cref{theorem: finite dual on algebraic K theory of ring} are inverse to each other and induce
   a natural equivalence of commutative ring spectra:
    \[
    \cK(A^\circ)\simeq  K(A).
    \]
    The equivalence is compatible with the Hattori--Stallings traces:
    \[
    \begin{tikzcd}
    \cK_0(A^\circ) \ar{rr}{\cong} \ar{dr}[swap]{\cotr}& & K_0(A)\ar{dl}{\tr}\\
   & A.
    \end{tikzcd}
    \]
\end{theorem}

\begin{proof}
    Given a finitely generated and projective $A$-module $M$, we show that the natural map $M\to M^{\circ *}$ is an isomorphism.
    Since $M$ is finitely generated, there exists $n$ such that there is a surjective $A$-linear homomorphism $A^{\oplus n}\rightarrow M$.
    As $M$ is projective, the previous map splits by a section $s\colon M\rightarrow A^{\oplus n}$. 
    We see that $M^{\circ *}$ inherits therefore a surjective map $(A^{\oplus n})^{\circ *}\cong (A^{\circ *})^{\oplus n}\to M^{\circ *}$ that splits with ${s^{\circ}}^*\colon M^{\circ *}\rightarrow (A^{\circ *})^{\oplus n}$. Therefore the natural map $M\rightarrow M^{\circ *}$ is a retract of the map $A^{\oplus n}\rightarrow (A^{\circ *})^{\oplus n}$, which is an isomorphism as $A$ is reflexive: 
    \[
    \begin{tikzcd}
    M^{\circ *} \ar[hook]{r}{s^{\circ *}} \ar[equals, bend left]{rr}& (A^{\circ *})^{\oplus n} \ar{r} & M^{\circ *}\\
    M \ar{u} \ar[hook]{r}{s} \ar[equals, bend right]{rr} & A^{\oplus n} \ar{u}{\cong} \ar{r} & M.\ar{u}
    \end{tikzcd}
    \]

Since $A$ is reflexive, its finite dual $A^\circ$ is coreflexive. 
    Therefore, we can argue similarly that given a finitely cogenerated and injective $A^\circ$-comodule $N$, the natural map $N\rightarrow N^{*\circ}$ is a retract of an isomorphism:
    \[
    \begin{tikzcd}
    N^{*\circ} \ar[hook]{r}\ar[equals, bend left]{rr}& (A^{\circ *\circ })^{\oplus n} \ar{r} & N^{*\circ}\\
    N \ar{u} \ar[hook]{r} \ar[equals, bend right]{rr} & (A^\circ)^{\oplus n} \ar{u}{\cong} \ar{r} & N.\ar{u}
    \end{tikzcd}
    \]
    Therefore  $(-)^*\colon \Inj_\fc(A^\circ)^\op\rightarrow \Proj_\fg(A)$ and $(-)^\circ\colon \Proj_\fg(A)\rightarrow \Inj_\fc(A^\circ)^\op$ are inverse of each others.
\end{proof}

\begin{corollary}\label{corollary: equivalence  K(C)=K(C^*)}
    Let $C$ be a $k$-coalgebra  such that its dual $C^*$ is a commutative Noetherian $k$-algebra.
    Then 
    the natural maps of \cref{thm coalg KT to alg KT} and \cref{corollary: inverse map K(C^*) to K(C)} are inverse to each other and induce a natural equivalence of commutative ring spectra:
    \[
    \cK(C)\simeq  K(C^*).
    \]
    The equivalence is compatible with the Hattori--Stallings traces:
    \[
    \begin{tikzcd}
    \cK_0(C) \ar{rr}{\cong} \ar{dr}[swap]{\cotr}& & K_0(C^*)\ar{dl}{\tr}\\
   & C^*.
    \end{tikzcd}
    \]
\end{corollary}

\begin{corollary}\label{corollary: K-theory of power series}
   There is an equivalence of commutative ring spectra:
    \[
    \cK(\dkx)\simeq K(k[[y]]).
    \]
    If $\mathrm{char}(k)=0$, there is an equivalence of commutative ring spectra:
    \[
    \cK(k[x])\simeq K(k[[y]]).
    \]
\end{corollary}

\begin{remark}\label{Remark: HH vs coHH on k[[y]] and k[x]}
While there is an equivalence between $K^c(k\langle X \rangle)$ and $K(k[[y]])$, the relationship between the coHochschild homology of the divided power coalgebra $k\langle X \rangle$ and the Hochschild homology of the power series ring $k[[y]]$ is more subtle.
Recall from \cite{Doi} or \cite[4.3]{cotrace} that $\coHH_n(C)$ is  $\mathrm{coTor}^n_{C\otimes C^\op}(C,C)$, where $C^\op$ denotes the coalgebra $C$ with reverse comultiplication. Consider the case where $\mathrm{char}(k)=0$, so that $k\langle X\rangle \cong k[x]$ (see \cref{example: k[x] is coreflexive}). It is shown in \cite[\S 2.3]{coHHComput} that:
    \[
    \coHH_n(k[x])= \begin{cases}
        k[x] & \text{if } n=0\\
        k[x]\otimes \langle x \rangle & \text{if }n=1\\
        0 & \text{if }n\geq 2,
    \end{cases}
    \]
    while the Hochschild homology of the ring $k[[y]]$ is given by:
    \[
    \HH_n(k[[y]])=\begin{cases}
        k[[y]] & \text{if }n=0\\
        k[[y]]\oplus W & \text{if }n=1\\
        \bigwedge^{n-1} W \oplus \bigwedge^{n} W & \text{if }n\geq 2,
    \end{cases}
    \]
where $W$ is the first Hochschild homology of $k((y))$ over $k(y)$, and in particular is an infinite dimensional $k((y))$-vector space, see \cite[2.5]{ayelet}.
In fact, $\coHH(k[x])$ more closely resembles the Hochschild homology of $k[[y]]$ seen as a pseudocompact algebra, see \cite[3.10]{IM25}.
\end{remark}

\section{Application to Swan theory of representations}
\label{sec: app to swan theory}

In the previous section, we showed that if $A$ is a $k$-algebra, then there is an equivalence $G^c(A^\circ)\simeq G^k(A)$ between the $G$-theories. 
In this section we consider the consequences of this equivalence in the case where $A$ is a group ring or monoid ring.

\begin{construction}
    Let $\Gamma$ be a monoid.
    Let $k\Gamma$ be the monoid ring with the multiplication in its $k$-algebra structure given by:
    \[
    \left( \sum_{g\in \Gamma} a_g g\right)\left( \sum_{g'\in \Gamma} b_{g'} g'\right)= \sum_{g\in \Gamma}\sum_{g'\in \Gamma}
     a_g b_{g'} (gg').\]
     It is endowed with a cocommutative $k$-coalgebra structure where 
    \begin{align*}
    \Delta\colon k\Gamma &\longrightarrow k\Gamma\otimes k\Gamma & \varepsilon\colon k\Gamma &\longrightarrow k \\
    g & \longmapsto g\otimes g & g & \longmapsto 1. 
    \end{align*}
    It is straightforward to check that $k\Gamma$ equipped with this multiplication and comultiplication is a cocommutative $k$-bialgebra. 
    
    In this section, we will consider the $G$-theory, $G^k(k\Gamma)$. When $\Gamma$ is a group, $G^k(k\Gamma)$ is referred to as the \emph{Swan theory} of $k\Gamma$ \cite{swan}. 
    If $\Gamma$ is a group, then $k\Gamma$ is a Hopf algebra with antipode $S(g)=g^{-1}$.
    \end{construction}

\begin{construction}\label{ConFin}
    Denote by $\Map(\Gamma, k)$ the vector space of all set functions $\Gamma\rightarrow k$.
    As there is an isomorphism $\Map(\Gamma, k)\cong k\Gamma^*$, then $\Map(\Gamma, k)$ inherits a commutative $k$-algebra structure, given by pointwise multiplication.
    Consider the $k$-linear morphism:
\[
\pi\colon \Map(\Gamma, k)\otimes \Map(\Gamma, k) \hookrightarrow \Map(\Gamma\times \Gamma, k)
\]
sending an element $f\otimes f'$ to a map $\Gamma\times \Gamma\rightarrow k$ defined by $(g,g')\mapsto f(g)f'(g')$.
If $\Gamma$ is finite, then $\pi$ is an isomorphism, and thus $\Map(\Gamma, k)$ inherits a $k$-coalgebra structure:
\[
\Delta\colon \Map(\Gamma, k)\longrightarrow  \Map(\Gamma\times \Gamma, k) \stackrel{\cong}\leftarrow \Map(\Gamma, k)\otimes \Map(\Gamma, k)
\]
where $\Delta$ assigns a function $f\colon \Gamma\rightarrow k$ to the function $\Delta f\colon \Gamma\times \Gamma\rightarrow k$ defined as $\Delta f(g, g')=f(gg')$, for all $g,g'\in \Gamma$.
\end{construction}

In general, if $\Gamma$ is not finite, then there is no $k$-coalgebra structure on $\Map(\Gamma, k)$, we need to instead consider the Sweedler dual $k\Gamma^\circ$ as we now explain.

\begin{construction}
For any monoid $\Gamma$, we can view $\Map(\Gamma, k)$ as a two-sided $k\Gamma$-module, where for all functions $f\colon \Gamma\to k$, and all $g,g',h\in \Gamma$, the two-sided action is determined by:
\[
(g\cdot f \cdot g')(h)\coloneqq f(g'hg).
\]
There is a $k$-algebra morphism
\[
\delta \colon \Map(\Gamma, k) \longrightarrow \Map(\Gamma\times \Gamma, k)
\]
that assigns a function $f\colon \Gamma\to k$ to a function $\Gamma\times \Gamma\rightarrow k$ defined by $(g,g')\mapsto f(gg')$. The morphism $\pi$ of \cref{ConFin} is in fact a $k$-algebra homomorphism in the sense that given functions $f_1,f_2,f'_1,f'_2\colon \Gamma\to k$, there is an equality:
\(
\pi(f_1f_1'\otimes f_2f_2')=\pi(f_1\otimes f_2)\pi(f_1'\otimes f_2').
\)
A classical result states that $\delta(f) \subseteq \pi\left(\Map(\Gamma, k)\otimes \Map(\Gamma, k)\right)$ if and only if the smallest $k\Gamma$-submodule of $\Map(\Gamma, k)$ containing $f$ is finite dimensional  (see \cite[Theorem 2.2.7]{abe}). 
Such a function $f\colon \Gamma\rightarrow k$ is called a \textit{representative function} of $\Gamma$. 
Since $\delta$ and $\pi$ are $k$-algebra homomorphisms, if $f,f'\colon \Gamma\to k$ are representative functions, then $\delta(ff')=\delta(f)\delta(f')$ and must be in the image of the $k$-subalgebra $\pi(\Map(\Gamma, k)\otimes \Map(\Gamma, k))$. Therefore $ff'$ must also be a representative function.
\end{construction}

\begin{definition}
    Define $R_k(\Gamma)$ to be the (commutative) subalgebra of $\Map(\Gamma, k)$ of representative functions.
    By construction, it is endowed with a $k$-coalgebra structure determined by $\delta$ above. We refer to $R_k(\Gamma)$ as the \emph{representative bialgebra of $\Gamma$}. 
    If $\Gamma$ is a group then $R_k(\Gamma)$ is a Hopf algebra. See \cite{abe} for more details on representative bialgebras.
\end{definition}

It is not hard to see that if $\Gamma$ is finite, then $R_k(\Gamma)\cong k\Gamma^*$ as bialgebras.

\begin{remark}
    If $\Gamma$ is only a semigroup, then $k\Gamma$ is a non-unital $k$-algebra, and $R_k(\Gamma)$ is a non-counital $k$-coalgebra.
    Given a $k$-algebra $A$, we can view it as a semigroup with its multiplication, yielding:
    \[
    A^\circ \cong R_k(A) \cap A^*.
    \]
    Thus $A^\circ$ is a subcoalgebra of $R_k(A)$.
    In particular, $R_k(\Gamma)= k\Gamma^\circ$. See also \cite[\S 2]{Taf72}. 
\end{remark}



\begin{example}
    Let $k$ be a finite field, and $\Gamma$ an infinite group that does not have any proper subgroup of finite index (e.g.\ a simple infinite group). Then $R_k(\Gamma)=k$, see \cite[2.3]{Taf72}. 
\end{example}

\begin{example}
    A locally finite representation of $\Gamma$ is a locally  finite left $k\Gamma$-module, and thus a right $R_k(\Gamma)$-comodule. Indeed, every locally finite 
     left $k\Gamma$-module $V$ with action $\rho\colon k\Gamma\otimes V\to V$ admits a natural $R_k(\Gamma)$-comodule structure:
\begin{align*}
    V & \longrightarrow V\otimes R_k(\Gamma)\\
    e_i & \longmapsto \sum_j e_j\otimes \rho_{ji}
\end{align*}
where $\{e_i\}_{i\in I}$ is a basis of $V$ and 
\[
\rho(g\otimes e_i)=\sum_{j}\rho_{ji}(g)e_j
\]
for all $g\in \Gamma$. The sum above is necessarily finite as $V$ is locally finite.
\end{example}

Often referred to as the algebraic $K$-theory of (finite dimensional) representations of $\Gamma$, recall that the Swan theory of $\Gamma$ is the algebraic $K$-theory $G^k(k\Gamma)$.
By \cref{corollary: equivalence G(A)=G(A^*)}, we see we can study Swan $K$-theory via the coalgebra $R_k(\Gamma)$.
We can also encode character theory co-algebraically.
Recall that the character of a finite dimensional representation $\rho\colon \Gamma\to \GL_n(k)$ of a group $\Gamma$ is a function $\chi(\rho)\colon \Gamma\to k$ where $\chi(\rho)(g)=\tr(\rho(g))$.
It defines a homomorphism $\chi\colon G^k_0(k\Gamma)\to \Map(\Gamma, k)$. Therefore \cref{corollary: equivalence G(A)=G(A^*)} induces the following.

\begin{corollary}\label{cor: Swan theory as coalgebras}
    Let $\Gamma$ be any group.
    There is an equivalence of algebraic $K$-theory spectra:
    \[
    \bK(R_k(\Gamma))\simeq G^k(k\Gamma),
    \]
    compatible with its character homomorphism in the sense that  the following diagram is commutative:
    \[
    \begin{tikzcd}
    \bK_0(R_k(\Gamma)) \ar{r}{\chi^c} \ar{d}[swap]{\cong} & \coHH_0(R_k(\Gamma)) \ar[hook]{d}\\
    G^k_0(k\Gamma) \ar{r}{\chi} & \Map(\Gamma, k).
    \end{tikzcd}
    \]
\end{corollary}

\begin{remark}
We can also now introduce an algebraic $K$-theory of locally finite representations of $\Gamma$ as $\cK(R_k(\Gamma))$. 
If $\Gamma$ is a Noetherian and abelian group, there are maps that fit in the following commutative diagram:
\[
\begin{tikzcd}
K_0(k\Gamma) \ar{r}\ar{d}[swap]{\tr} & \cK_0(R_k(\Gamma))\ar{d}{\cotr}\\
k\Gamma \ar{r} & R_k(\Gamma)^*.
\end{tikzcd}
\]
\end{remark}

\bibliographystyle{alpha}
\bibliography{bib}

\end{document}